\documentclass[10pt]{article}
\usepackage[top=1in, bottom=1in, left=1in, right=1in]{geometry}
\usepackage{srcltx,graphicx}
\usepackage{amsmath, amssymb, amsthm}
\usepackage{mathtools}
\usepackage{color, xcolor, colortbl}
\usepackage{subfig, overpic}
\usepackage{hyperref}
\usepackage{epstopdf}
\usepackage{algorithm}
\usepackage{algpseudocode}
\algnewcommand{\LineComment}[1]{\Statex \hspace{-5mm}\(\triangleright\) #1}
\usepackage[capitalize,nameinlink]{cleveref} \usepackage{caption}
\usepackage{epstopdf}
\usepackage{tikz}
\usetikzlibrary{math} 

\usepackage{stackengine}
\usepackage{scalerel}

\graphicspath{{images/}}
\numberwithin{equation}{section}

\newcommand\bbR{\mathbb{R}}
\newcommand\bbZ{\mathbb{Z}}
\newcommand\bbN{\mathbb{N}}
\newcommand\bbS{\mathbb{S}}

\newcommand\dd{\,\mathrm{d}} \newcommand\diag{\mathrm{diag}\,}
\newcommand\sps[1]{^{(#1)}}
\newcommand\wl[2]{\sps{#1}_{#2}}

\newcommand\cM{\mathcal{M}}
\newcommand\cN{\mathcal{N}}

\newcommand\cU{\mathcal{U}}

\newcommand\K{n_{\mathrm{cnn}}}

\newcommand\fw{\psi}
\newcommand\fs{\varphi}
\newcommand\cw{d}
\newcommand\cs{v}

\newcommand\T{{\mathsf{T}}}

\newcommand\Gv{{G_\eta}}

\newcommand\Gmu{{\Gv}}

\newcommand\sffont[1]{{\sf{#1}}}

\newcommand\id{{\sffont{id}}}

\newcommand\ConvNet{{\sffont{ConvNet} }}

\newcommand\Convone{{\sffont{Conv1d }}}
\newcommand\Convtwo{{\sffont{Conv2d }}}

\newcommand\FWT{\sffont{FWT}}
\newcommand\IWT{\sffont{IWT}}

\newcommand\Nparams{{N_{\mathrm{params}}}}
\newcommand\Nsamp{{N_{\mathrm{samples}}}}
\newcommand\trainerror{{\epsilon_{\mathrm{train}}}}
\newcommand\testerror{{\epsilon_{\mathrm{test}}}}
\newcommand\operror{{\epsilon_{\mathrm{op}}}}
\newcommand\NN{\mathrm{NN}}

\title{
  Meta-learning Pseudo-differential Operators \\with Deep Neural Networks
}

\date{}
\author{
  Jordi Feliu-Fab{\`a}  \thanks{ICME, Stanford University, Stanford, CA 94305. 
    Email: {\tt jfeliu@stanford.edu}},~~
  Yuwei Fan  \thanks{Department of Mathematics, Stanford University, Stanford, CA 94305. 
    Email: {\tt ywfan@stanford.edu}},~~
  Lexing Ying  \thanks{Department of Mathematics and ICME, Stanford University, Stanford, CA 94305. 
    Email: {\tt lexing@stanford.edu}}
}

\newcommand{\RRev}[1]{#1}

\newcommand{\grad}{\nabla}
\newcommand{\lapl}{\Delta}

\newcommand{\bi}{\begin{itemize}}
\newcommand{\ei}{\end{itemize}}

\begin{document}
\maketitle
\begin{abstract}
  This paper introduces a meta-learning approach for parameterized pseudo-differential operators
  with deep neural networks. With the help of the nonstandard wavelet form, the pseudo-differential
  operators can be approximated in a compressed form with a collection of vectors. The nonlinear map
  from the parameter to this collection of vectors and the wavelet transform are learned together
  from a small number of matrix-vector multiplications of the pseudo-differential operator.
  Numerical results for Green's functions of elliptic partial differential equations and the
  radiative transfer equations demonstrate the efficiency and accuracy of the proposed approach.

                    \end{abstract}


\section{Introduction}\label{sec:intro}
Many physical models for scientific and engineering applications can be written in a general form
\begin{equation}
  L_{\eta} u(x)=f(x), \quad x\in \Omega \subset \bbR^d
  \label{eq:PDE}
\end{equation}
for a domain $\Omega$ with appropriate boundary conditions, where $L_{\eta}$ is often a partial
differential or integral operator parameterized by a parameter function $\eta(x)$. Solving for $u(x)$
for a given $f(x)$ amounts to representing the inverse operator (sometimes also known as the Green's
function) $G_\eta = L_{\eta}^{-1}$ either explicitly or implicitly via an efficient
algorithm. Representing $G_\eta$, even if implicitly, can be computationally challenging,
especially for multidimensional problems. The past few decades have witnessed steady progresses in
developing efficient algorithms for this. 

\paragraph{Problem statement.}
This paper is concerned with a more ambitious task: representing the nonlinear map from $\eta$ to
$G_\eta$
\begin{equation}
  \cM: \eta \rightarrow G_\eta = L_{\eta}^{-1},
  \label{eq:Mmap}
\end{equation}
when the operator $G_\eta$ is a pseudo-differential operator (PDO) \cite{wong2014introduction}.
Although $L_{\eta}$ and $G_\eta$ can be linear operators, this map $\cM$ from $\eta$ to $G_\eta$ is
highly nonlinear.

\paragraph{Background.}
In the recent years, deep learning has become the most versatile and effective tool in artificial
intelligence and machine learning, witnessed by impressive achievements in computer vision
\cite{Krizhevsky2012,Zeiler2014,He2016DeepRL}, speech and natural language processing
\cite{Hinton2012,Socher12,Sarikaya14, vaswani2017attention,devlin2018bert}, drug discovery
\cite{MaSheridan2015} or game playing \cite{Silver2016,evans2018novo,vinyals2017starcraft}. Recent
reviews on deep learning and its impacts on other fields can be found in for example
\cite{leCunn2015,SCHMIDHUBER2015}. At the center of deep learning, the model of deep neural networks
(NNs) provides a flexible framework for approximating high-dimensional functions, while allowing for
efficient training and good generalization properties in practice \cite{Livni14,Poggio16}.

More recently, several groups have started applying NNs to partial differential equations (PDEs) and
integral equations (IEs) arising from physical systems. In one direction, the NN model has been used
to approximate solutions of high-dimensional PDEs
\cite{lagiaris98, sirigano17, DeepRitz2017, Karni2017, berg2017unified, carleo2017solving,
Weinan2017, han2018solving, khoo2019committor, Liu2019}.
In a somewhat orthogonal direction, the NNs have been utilized to approximate the high-dimensional
parameter-to-solution of various PDEs and IEs \cite{khoo2017solving, han2017deep, fan2018mnn,
fan2018mnn2, bcrnet, switchnet, han2018solving, Araya-Polo2018, li2019variational,fan2019eit}.

Another topic from machine learning that is particularly relevant to this work is {\em
meta-learning} or {\em learning-to-learn} \cite{meta87,Bengio1990,lowshot2017,Finn2017,wang16}.  A
meta-learning system learns to produce learning models for new tasks and scenarios from their
metadata with zero or minimum amount of new data, by leveraging the common structure among different
tasks. Due to the low requirements on new data points, meta-learning has gained a lot of attention
in recent years in applications such as vision and reinforcement learning.

\paragraph{Main idea.}
Following these recent advances in applying NNs to physical models, this paper takes a deep learning
approach for representing the map in \cref{eq:Mmap}. The most straightforward solution would be to
take a supervised learning approach, i.e., trying to learn the map $\mathcal{M}: \eta\rightarrow
G_\eta$ from a large set of training data $\{(\eta_i,G_{\eta_i})\}_i$. However, since it is often
difficult or even impossible to compute and store $G_{\eta_i}$ due to the enormous discretization
size, this straightforward supervised learning approach is not practical for \cref{eq:Mmap}.

Without explicit access to $G_\eta$, we take a meta-learning approach, i.e., learning to produce,
for each new $\eta$, an NN approximation to $G_{\eta}$. To do this, we are faced with two key
difficulties.
\begin{itemize}
\item How should we represent the output $G_\eta$ for an arbitrary input $\eta$?
\item How should we represent the training data?
\end{itemize}

To address the first question, $G_\eta$ should be represented in a compressed form. For
pseudo-differential operators, several compressed representations exist, including hierarchical
matrices \cite{hackbusch1999sparse,hackbusch2001introduction,Hackbusch2000H2}, discrete symbol
calculus \cite{Symb2011}, etc.
In this paper, we choose to represent $G_\eta$ with the nonstandard wavelet form introduced in
\cite{bcr}. The main advantage of the nonstandard wavelet form is that the nonzero entries of this
compressed representation are simply organized into a small number of vectors. More precisely,
\[
G_\eta \approx W \mathcal{S}[C_\eta] W^\T,
\]
where $W$ is a redundant form of a wavelet transform, $C_\eta$ stands for the collection of vectors
that contain the nonzero entries of the compressed form, and $\mathcal{S}$ is a certain operator
that generates a sparse matrix from the vector collection $C_\eta$. Compared to \cite{bcr}, a key
difference is that the current approach allows for $W$ to be fine-tuned for the map $\mathcal{M}$.

To address the second question, instead of explicitly representing $G_{\eta_i}$, the training data
consists of samples of the form
\[
(\eta_i, \{f_{ij}, u_{ij}\}),
\]
where $u_{ij} = G_{\eta_i} f_{ij}$. For a fixed $\eta_i$, such data can be obtained by solving the
equation $L_{\eta_i} u_{ij} = f_{ij}$ for each $f_{ij}$, possibly with a fast algorithm.

Putting these two pieces together, the meta-learning approach of this paper learns two following key
objects from the training data of form $\{ (\eta_i, \{f_{ij}, u_{ij}\}_j) \}_i$:
\begin{itemize}
\item a map from $\eta$ to the vector collection $C_\eta$,
\item the $\mathcal{M}$-dependent wavelet transform $W$.
\end{itemize}
Once trained, for a given test input $\eta$ the architecture calculates $C_\eta$ and returns a
linear NN that implements $W \mathcal{S}[C_\eta] W^\T \approx G_\eta$.

\paragraph{Organization.}
The rest of this paper is organized as follows. \cref{sec:Pre} briefly reviews the nonstandard
wavelet form, used for representing $G_{\eta}$. In \cref{sec:nn}, the NN architecture of the
meta-learning approach is discussed in detail. \Cref{sec:sch} applies the proposed NN to the
Green's function of elliptic PDEs, in both the Schr{\"o}dinger form and the divergence form.
The application to the radiative transfer equation is presented in \cref{sec:rte}.

\section{Nonstandard wavelet form}\label{sec:Pre}

This section summarizes the nonstandard wavelet form proposed in \cite{bcr}. To make things
concrete, compactly supported orthonormal Daubechies wavelets \cite{daubechies1988orthonormal} are
used as the basis functions as an example.

\subsection{Wavelet transform}
In the one-dimensional multiresolution analysis, one starts by defining a scaling function
$\varphi(x)$ that generates, through dyadic translations and dilations, a family of functions
\begin{equation}\label{eq:translate_phi}
  \varphi_{k}^{(\ell)}(x)=2^{\ell/2}\varphi(2^{\ell}x-k),\quad
  \ell=0,1,2,\dots,\quad k\in\bbZ.
\end{equation}
For each scale $\ell$, the functions $\{\varphi_{k}^{(\ell)}\}$ form a Ritz basis for a space
$V_{\ell}$, which satisfies a nested relationship $V_{\ell}\subset V_{\ell+1}$. This nested property
of $\{V_{\ell}\}$ implies the following \emph{dilation relation} of the scaling function
\begin{equation}\label{eq:recursive_phi}
  \varphi(x) = \sqrt{2} \sum_{i\in\bbZ} h_i\varphi(2x-i).
\end{equation}
For the Daubechies' wavelets \cite{daubechies1988orthonormal}, the scaling function $\varphi(x)$ has
a compact support $[0, 2p-1]$ for a given positive integer $p$ and therefore the coefficients
$\{h_i\}$ are only nonzero for $i=0,\dots,2p-1$. The scaling function also satisfies the orthonormal
condition
\begin{equation}\label{eq:normalized}
  \int_{\bbR}\varphi(x-a)\varphi(x-b)\dd x = \delta_{a,b}, \quad \forall a,b\in\bbZ,
\end{equation}
which leads to an orthonormal condition for the coefficients $\{h_i\}$
\begin{equation}\label{eq:orthogonality}
  \sum_{i\in\bbZ} h_i^2=1,\quad 
  \sum_{i\in\bbZ} h_ih_{i+2m}=0, \quad m\in\bbZ \backslash \{0\}.
\end{equation}

Given the scaling function $\varphi(x)$, another important component of the multiresolution analysis
is the wavelet function $\psi(x)$, defined by
\begin{equation}\label{eq:recursive_psi}
  \psi(x) = \sqrt{2} \sum_{i\in\bbZ} g_i\varphi(2x-i),
\end{equation}
where $g_i=(-1)^{1-i} h_{1-i}$ for $i\in\bbZ$.  A simple calculation shows that the support of
$\psi(x)$ is $[-p+1,p]$ and $\{g_i\}$ is nonzero only for $i=-2p+2,\ldots,1$, based on the support
of the $\varphi$ and the nonzero entries pattern of $\{h_i\}$.  The {\em Daubechies wavelets} are
then defined as
\begin{equation} \label{eq:cell_psi}
  \psi\wl{\ell}{k}(x)    = 2^{\ell/2}\psi\left( 2^{\ell}x-k \right),
  \quad \ell=0,1,2,\dots,\quad k\in\bbZ.
\end{equation}

For a function $v(x)\in L^2(\bbR)$, its scaling and wavelet coefficients $\cw_k\sps{\ell}$ and
$\cs_k\sps{\ell}$ are defined as the inner product with the scaling functions and the wavelets 
\begin{equation}
  \cs_k\sps{\ell} := \int v(x)\fs_k\sps{\ell}(x)dx, \quad
  \cw_k\sps{\ell} := \int v(x)\fw_k\sps{\ell}(x)dx.
\end{equation}
Using the recursive relationships of the scaling function \cref{eq:recursive_phi} and the wavelet
function \cref{eq:recursive_psi}, one obtains a recursive relationship of the scaling and wavelet
coefficients
\begin{equation}
  \label{eq:coeff}  
  \cs_k\sps{\ell} = \sum_{i\in \bbZ}h_i \cs_{2k+{i}}^{(\ell+{1})},\quad
  \cw_k\sps{\ell} = \sum_{i\in \bbZ}g_i \cs_{2k+{i}}^{(\ell+{1})}.
\end{equation}
By defining $\cs\sps{\ell}=\left(\cs\sps{\ell}_k\right)_{k\in\bbZ}$ and
$\cw\sps{\ell}=\left(\cw\sps{\ell}_k\right)_{k\in\bbZ}$, \cref{eq:coeff} can be written in a matrix
form 
\begin{equation}\label{eq:v2vd}
  \cs\sps{\ell}=\left(W_s\sps{\ell}\right)^\T\cs\sps{\ell+1},\quad
  \cw\sps{\ell}=\left(W_w\sps{\ell}\right)^\T\cs\sps{\ell+1},
\end{equation}
where the operators $W_s\sps{\ell}$ and $W_w\sps{\ell}:\ell^2(\bbZ)\to\ell^2(\bbZ)$ are banded with
a bandwidth $2p$ due to the support of $\{h_i\}$ and $\{g_i\}$. By introducing the orthogonal
operator $W\sps{\ell} = \left(W_w\sps{\ell} \quad W_s\sps{\ell}\right)$, \cref{eq:v2vd} can be
rewritten as
\begin{equation} \label{eq:FWT}
    \begin{pmatrix}
        \cw\sps{\ell} \\ 
        \cs\sps{\ell}
    \end{pmatrix}
    = \left(W\sps{\ell}\right)^\T \cs^{(\ell+1)},
    \quad
    \cs\sps{\ell+1}= W\sps{\ell} 
    \begin{pmatrix}
        \cw\sps{\ell}\\
        \cs^{(\ell)}
    \end{pmatrix}.
\end{equation}
The procedure for computing the wavelet and scaling coefficients can be illustrated in the following
diagram
\begin{equation}\label{eq:procedure}
  \newcommand\lra{\longrightarrow}
  \setlength{\arraycolsep}{2pt}
  \begin{array}{ccccccccccccccc}
    \cdots & \lra & \cs\sps{\ell}& \lra     & \cs\sps{\ell-1} & \lra     & \cs\sps{\ell-2} & \lra
    & \cdots   & \lra     & \cs\sps{2} & \lra       & \cs\sps{1} & \lra      & \cs\sps{0} \\
    & \searrow &      & \searrow &          & \searrow &            & \searrow
    &          & \searrow &          & \searrow   &          & \searrow \\
    & &  \cw\sps{\ell}    &          & \cw\sps{\ell-1} &          & \cw\sps{\ell-2} &
    & \cdots   &          & \cw\sps{2} &            & \cw\sps{1} &           & \cw\sps{0}
  \end{array}.
\end{equation}

The discussion until now is concerned with the wavelets on $\bbR$. It is straightforward to extend
it the functions defined on a finite domain with periodic boundary condition. If the function $v(x)$
is periodic on a finite domain, for instance, $[0,1]$, then the only modification is that all the
shifts and scaling in the $x$ variable are done modulus the integer.
When working with periodic functions, the procedure in \cref{eq:procedure} usually stops at a coarse
level $L_0=O(\log_2(p))$ before the wavelet and scaling functions start to overlap itself.
\begin{equation}
    \newcommand\lra{\longrightarrow}
    \setlength{\arraycolsep}{2pt}
    \begin{array}{ccccccccccccc}
         & &\cdots& \lra & \cs\sps{\ell+1}& \lra     & \cs\sps{\ell} & \lra     & \cs\sps{\ell-1} & \lra
            & \cdots   & \lra        & \cs\sps{L_0} \\
        && &\searrow &      & \searrow &          & \searrow &            & \searrow
            &          & \searrow \\
        & &&&  \cw\sps{\ell+1}    &          & \cw\sps{\ell} &          & \cw\sps{\ell-1} &
            & \cdots   &        & \cw\sps{L_0}
    \end{array}.
\end{equation}

\subsection{Nonstandard wavelet form for integral operator}

Let $A$ be an integral operator with kernel $a(x,y)$, applied to periodic functions defined on
$[0,1]$, i.e.,
\begin{equation}\label{eq:Av_operator}
  u = A v,\quad\text{equivalently}\quad  u(x) = \int a(x,y) v(y)\dd y.
\end{equation}
Denote by $A\sps{L}=\left( A_{k_1,k_2}\sps{L} \right) \in\bbR^{2^L\times 2^L}$ the Galerkin
projection of $A$ to the space $V_{L}$, for a sufficiently deep level $L$, i.e.
\[
A_{k_1,k_2}\sps{L} = \int\int \varphi_{k_1}\sps{L}(x) a(x,y) \varphi_{k_2}\sps{L}(y) \dd x\dd y.
\]
The nonstandard form described in \cite{bcr} is a remarkably efficient way to compress the matrix
$A\sps{L}$.

The main step for the nonstandard form is to treat $A\sps{L}$ as an image and use the 2D
multiresolution analysis 
\begin{equation}\label{eq:Int_D2A}
  \begin{aligned}
    D\wl{\ell}{1,k_1,k_2}&:=\iint \psi\wl{\ell}{k_1}(x)a(x,y)\psi\wl{\ell}{k_2}(y)\dd x\dd y,
    &
    D\wl{\ell}{2,k_1,k_2}&:=\iint \psi\wl{\ell}{k_1}(x)a(x,y)\varphi\wl{\ell}{k_2}(y)\dd x\dd y,\\
    D\wl{\ell}{3,k_1,k_2}&:=\iint \varphi\wl{\ell}{k_1}(x)a(x,y)\psi\wl{\ell}{k_2}(y)\dd x\dd y,
    &
    A\wl{\ell}{k_1,k_2}  &:=\iint \varphi\wl{\ell}{k_1}(x)a(x,y)\varphi\wl{\ell}{k_2}(y)\dd x\dd y,
  \end{aligned}
\end{equation}
for $\ell=L_0, \dots, L-1$, and $k_1,k_2 = 0,\dots,2^\ell-1$.  For convenience, these coefficients
are organized into the matrix form as
\begin{equation}\label{eq:AD}
  A\sps{\ell}=(A\sps{\ell}_{k_1,k_2})_{k_1,k_2=0,\dots,2^\ell-1},\quad\quad
  D_j\sps{\ell}=(D\sps{\ell}_{j,k_1,k_2})_{k_1,k_2=0,\dots,2^\ell-1}, j=1,2,3.
\end{equation}
In this setting, a similar recursive relation to \cref{eq:FWT} can be obtained
\begin{equation} \label{eq:FWT_2D}
    \begin{pmatrix}
        D_1\sps{\ell}   & D_2\sps{\ell}\\
        D_3\sps{\ell}   & A\sps{\ell}
    \end{pmatrix}
    = (W\sps{\ell})^\T A\sps{\ell+1} W\sps{\ell},\quad \ell=L_0,\dots,L-1.
\end{equation}
If $A$ is a Calderon-Zygmund operator, the entries of the matrices $D_j\sps{\ell}$ with $j=1,2,3$
decay rapidly away from the diagonal. For a prescribed relative accuracy $\epsilon$, each matrix
$D_j\sps{\ell}$ can be approximated by a band matrix by truncating at a band of width
$O(\log(1/\epsilon))$. Since the bandwidth is independent of the specific choices of $\ell$, $j$, or
the mesh size $N=2^L$, the nonstandard form of $A$ stores only $O(N)$ nonzero entries. The readers
are referred to \cite{bcr} for more details. With a slight abuse of notation, the matrices
$D_j\sps{\ell}$ are assumed to be pre-truncated in what follows.

\begin{figure}[ht]
  \centering
  \includegraphics[width=0.95\textwidth]{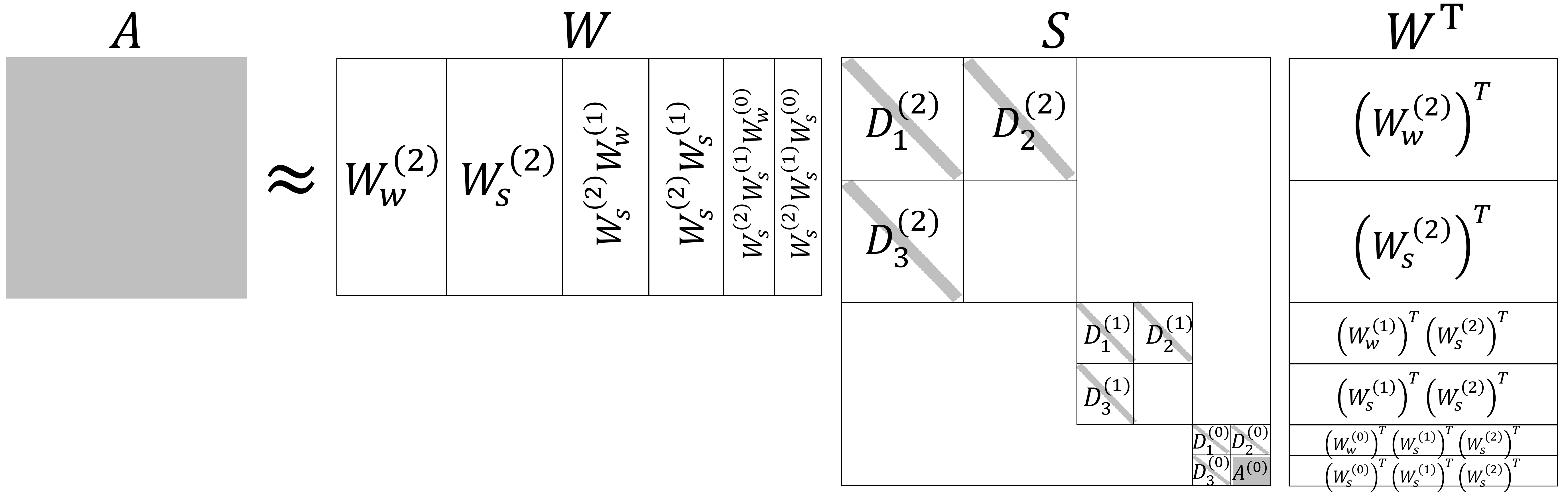}
  \caption{\label{fig:WHW}Illustration of $A=W S W^\T$ with $L_0=0$ and $L=3$. The sparsity
    pattern of $S$ is marked in gray.}
\end{figure}

One can assemble all the matrices $D_j^{\ell}$ and $A\sps{L_0}$ together, by defining the matrix
$S\sps{L}$ in a recursive way as
\begin{equation}
  S\sps{L_0} = \begin{pmatrix}
    D_1\sps{L_0} & D_2\sps{L_0} \\
    D_3\sps{L_0} & A\sps{L_0}
  \end{pmatrix},\quad
  S\sps{\ell+1} = \begin{pmatrix}
    D_1\sps{\ell+1} & D_2\sps{\ell+1} & 0 \\
    D_3\sps{\ell+1} &  0 & 0\\
    0 & 0 & S\sps{\ell}
  \end{pmatrix},\quad \ell=L_0,\dots,L-1.
\end{equation}
The matrix $S := S\sps{L}$ is the nonstandard form of the matrix $A=A\sps{L}$ satisfying
\begin{equation}\label{eq:WHW}
  A = W S W^\T.
\end{equation}
Here $W$ is the {\em extended} wavelet transform matrix, defined in the recursive form as
\begin{equation}
  T\sps{L_0} = W\sps{L_0},\quad 
  T\sps{\ell+1} = \left( W\sps{\ell+1} \quad W_s\sps{\ell+1}T\sps{L_0} \right),
  \quad \ell=L_0,\dots,L-1,
  \quad W := T\sps{L}.
\end{equation}
\Cref{fig:WHW} illustrates the matrices $W$ and $S$ along with the formulation \cref{eq:WHW}.

To clarify the notations, we denote by $W_s\sps{\ell}$ and $W_w\sps{\ell}$ the transform matrices
defined in \cref{eq:v2vd} for the scaling and wavelet parts on level $\ell$, respectively.
$W\sps{\ell}=(W_w\sps{\ell}\quad W_s\sps{\ell})$ is the wavelet transform matrix at the level $\ell$.

\subsection{Matrix-vector multiplication in the nonstandard form}\label{sec:matvec}

With a Galerkin discretization of \cref{eq:Av_operator} at level $L$, the matrix-vector
multiplication takes the form
\begin{equation}\label{eq:Av}
  u\sps{L}=A\sps{L}v\sps{L}.
\end{equation}
The nonstandard form allows for accelerating the evaluation of \cref{eq:Av}. Using the nonstandard
form $A\sps{L} = W S W^\T$ obtained above, the matrix-vector multiplication 
\begin{equation}
  u\sps{L}=W S W^\T v\sps{L}
\end{equation}
can be split into four steps:
\begin{enumerate}
\item $A\sps{L}\to S$: generate the nonstandard form $S$ from the matrix $A\sps{L}$ or the
  kernel $a(x,y)$;
\item $v\sps{L}\to \hat{v}:=W^\T v\sps{L}$: apply (forward) wavelet transform on $v\sps{L}$ to
  get $\hat{v}$;
\item $\hat{u} := S\hat{v}$: evaluate the matrix-vector multiplication in the nonstandard form;
\item $\hat{u}\to u\sps{L}:=W\hat{u}$: apply inverse wavelet transform on $\hat{u}$ to obtain
  $u\sps{L}$.
\end{enumerate}
The first step is computed using \cref{eq:FWT_2D} if the matrix $A\sps{L}$ is given. The second step
follows \cref{eq:FWT}. The third step can be written as
\begin{equation}
    \begin{pmatrix}
        w\sps{\ell}\\
        s\sps{\ell}
    \end{pmatrix}
    = \begin{pmatrix}
        D_1\sps{\ell}   & D_2\sps{\ell}\\
        D_3\sps{\ell}   & D_4\sps{\ell}
    \end{pmatrix}
    \begin{pmatrix}
        d\sps{\ell}\\
        v\sps{\ell}
    \end{pmatrix},
\end{equation}
where $D_4\sps{\ell}=0$ for $\ell=L_0+1,\dots,L-1$ and $D_4\sps{L_0}=A\sps{L_0}$.  The fourth step
is essentially an inverse wavelet transform, implemented as
\begin{equation}
    u\sps{L_0}=0,\quad 
    u\sps{\ell+1}=W\sps{\ell}
    \begin{pmatrix}
        w\sps{\ell}\\
        s\sps{\ell}+u\sps{\ell}
    \end{pmatrix},
    \quad \ell=L_0,\dots,L-1.
\end{equation}
A step-by-step description of these four steps are summarized in \cref{alg:BCR}.

\begin{algorithm}[htb]
  \begin{small}
    \begin{center}
      \begin{minipage}{0.5\textwidth}
        \begin{algorithmic}[1]
          \Require $A\sps{L}=A$, $v\sps{L}=v$, $0\leq L_0<L$;
          \Ensure $u=u\sps{L}$;
          \vspace{1mm}
          \LineComment{ {\bf Step 1:} generate nonstandard form of $A\sps{L}$}
          \For {$\ell$ from $L-1$ to $L_0$ by $-1$}\vspace{2mm}
          \State  
          $\begin{pmatrix}
            D_1\sps{\ell} &  D_2\sps{\ell}\\
            D_3\sps{\ell} & A\sps{\ell}
          \end{pmatrix}  = (W\sps{\ell})^\T A\sps{\ell+1}W\sps{\ell}$
          \EndFor
          \vspace{2mm}
          \LineComment{ {\bf Step 2:} forward wavelet transform on $v\sps{L}$}
          \For {$\ell$ from $L-1$ to $L_0$ by $-1$}\vspace{2mm}
          \State 
          $\begin{pmatrix}
            d\sps{\ell}\\
            v\sps{\ell}
          \end{pmatrix}=(W\sps{\ell})^\T v\sps{\ell+1}$ 
          \EndFor
          \algstore{break1}
        \end{algorithmic}
      \end{minipage}
      \begin{minipage}{0.48\textwidth}
        \begin{algorithmic}[1]
          \algrestore{break1}
          \LineComment{ {\bf Step 3:} matrix-vector multiplication}
          \State $D_4\sps{\ell}=0$ for $\ell = L_0+1,\cdots, L-1$ and $D_4\sps{L_0}=A\sps{L_0}$
          \For {$\ell$ from $L_0$ to $L-1$}\vspace{1mm}
          \State 
          $\begin{pmatrix}
            w\sps{\ell}\\
            s\sps{\ell}
          \end{pmatrix}=\begin{pmatrix}
          D_1\sps{\ell} & D_2\sps{\ell}\\
          D_3\sps{\ell} &  D_4\sps{\ell}
          \end{pmatrix}\begin{pmatrix}
            \cw\sps{\ell}\\
            \cs\sps{\ell}
          \end{pmatrix}$
          \EndFor  
          \vspace{2mm}
          \LineComment{ {\bf Step 4:} Inverse wavelet transform}
          \State $u\sps{L_0}=0$
          \For {$\ell$ from $L_0$ to $L-1$}\vspace{2mm}
          \State 
          $u\sps{\ell+1}= W\sps{\ell} \begin{pmatrix}
            w\sps{\ell}\\
            s\sps{\ell}+u\sps{\ell}
          \end{pmatrix}$
          \EndFor
          \State $\mathsf{return}$ $u = u\sps{L}$
        \end{algorithmic}
        
      \end{minipage}
    \end{center}
  \end{small}
  \caption{Compute $u = Av$ using the nonstandard form of the wavelet transform}
  \label{alg:BCR}
\end{algorithm}

\subsection{The multidimensional case}\label{sec:matvec_2d}
The matrix-vector multiplication in the nonstandard form can be easily extended to the
multidimensional case with the help of multidimensional orthogonal wavelets (see \cite{Mallat2009}
for more details). For instance, in the two-dimensional setting, one defines at each scale $\ell$
three different types of wavelets of the form
\begin{equation}
  \fw_{1,k}\sps{\ell}(x,y) = \fs_{k_1}\sps{\ell}(x)\fw_{k_2}\sps{\ell}(y), \quad
  \fw_{2,k}\sps{\ell}(x,y) = \fw_{k_1}\sps{\ell}(x)\fs_{k_2}\sps{\ell}(y), \quad
  \fw_{3,k}\sps{\ell}(x,y) = \fw_{k_1}\sps{\ell}(x)\fw_{k_2}\sps{\ell}(y),
\end{equation}
with $k=(k_1,k_2)\in\bbZ^2$. Using these three types of wavelets, the transform matrix at each scale
$\ell$ used in \cref{eq:FWT} is redefined to be $W\sps{\ell} = \left(W_{w,1}\sps{\ell} \quad
W_{w,2}\sps{\ell} \quad W_{w,3}\sps{\ell} \quad W_s\sps{\ell}\right)$. The 2D analog of
\cref{eq:FWT} contains three types of wavelet coefficients
\begin{equation} \label{eq:FWT_2d}
    \begin{pmatrix}
        \cw_1\sps{\ell} \\ 
        \cw_2\sps{\ell} \\
        \cw_3\sps{\ell} \\
        \cs\sps{\ell}
    \end{pmatrix}
    = \left(W\sps{\ell}\right)^\T \cs^{(\ell+1)},
    \quad
    \cs\sps{\ell+1}= W\sps{\ell} 
    \begin{pmatrix}
        \cw_1\sps{\ell}\\
        \cw_2\sps{\ell}\\
        \cw_3\sps{\ell}\\
        \cs^{(\ell)}
    \end{pmatrix}.
\end{equation}
Similarly, the recursive relation \cref{eq:FWT_2D} can be extended as well
\begin{equation}\label{eq:FWT_4D}
  \begin{pmatrix}
    D_1\sps{\ell}   & D_2\sps{\ell}   & D_3\sps{\ell} & D_4\sps{\ell}\\
    D_5\sps{\ell}   & D_6\sps{\ell}   & D_7\sps{\ell} & D_8\sps{\ell}\\
    D_9\sps{\ell}   & D_{10}\sps{\ell} & D_{11}\sps{\ell} & D_{12}\sps{\ell}\\
    D_{13}\sps{\ell} & D_{14}\sps{\ell}  & D_{15}\sps{\ell} & A\sps{\ell}
  \end{pmatrix}
  = (W\sps{\ell})^T A\sps{\ell+1} W\sps{\ell},
  \quad \ell=L_0,\dots,L-1,
\end{equation}
where $D_j\sps{\ell}$, $j=1,\dots,15$ are all sparse matrices with only $O(4^\ell)$ non-negligible
entries in each. The matrix-vector multiplication follows the steps of \cref{alg:BCR}, with these
necessary changes.

\section{Meta-learning approach}\label{sec:nn}

The plan is to apply the nonstandard form to the operator $G_\eta$ in \cref{eq:Mmap}
\begin{equation}\label{eq:Gf}
  u = G_\eta f,\quad 
  u(x) = \int g_{\eta}(x,y)f(y)\dd y.
\end{equation}
With a slight abuse of notations, the same letters are used to denote the discretizations. The
discrete version of \cref{eq:Gf} takes the form
\begin{equation}
  u = G_{\eta} f,\quad u, f, \eta\in\bbR^{N} \text{ and } G_{\eta} \in \bbR^{N\times N},
\end{equation}
with $N=2^L$. The main goal of this paper is to construct a neural network to learn the map $\eta\to
G_{\eta}$. 

Following \cref{eq:WHW} and applying the wavelet transform to the matrix $G_\eta$ leads to
\begin{equation}\label{eq:WHW_nn}
  G_{\eta}\approx W S_{\eta} W^\T,
\end{equation}
where $W$ is the extended wavelet transform matrix, independent of the parameter $\eta$. Since each
block of matrix $S_{\eta}$ is a band matrix, the nonzero entries of each block can be represented by
a set of vectors. Let us define $C_{\eta}\sps{\ell}$, of size $2^\ell\times n_c$, to be the
collection of these vectors of $S_{\eta}$ at level $\ell$, with $n_c$ dependent on the bandwidth and
$\ell$. By introducing the collection of vectors $C_{\eta} :=
\{C_{\eta}\sps{\ell}\}_{\ell=L_0,\dots,L-1}$, $S_{\eta}$ is uniquely determined by $C_{\eta}$, i.e.,
\[
S_{\eta} \equiv \mathcal{S}[C_{\eta}]
\]
for a fixed embedding operator $\mathcal{S}$ determined by the sparsity pattern of $S_{\eta}$.

Given a set of training samples of the form
\begin{equation}\label{eq:data}
  (\eta_i, \{f_{ij}, u_{ij}\}),
\end{equation}
where $u_{ij} = G_{\eta_i} f_{ij}$ can be obtained by solving $L_{\eta_i}u_{ij}=f_{ij}$ with right
hand side $f_{ij}$, the meta-learning approach first learns both the map $\eta\to C_\eta$ and the
wavelet transform matrix $W$. Once they are ready, given any new $\eta$, $G_{\eta}$ can be
approximated by evaluating the map $\eta\to C_{\eta}$ and representing \eqref{eq:WHW_nn} in an NN
form.

\subsection{Neural network architecture}

Using the factorization of $G_{\eta}$ in \cref{eq:WHW_nn}, one can factorize $u_{ij} =
G_{\eta_i}f_{ij}$ as
\begin{equation}\label{eq:WHWf}
  u_{ij} \approx W S_{\eta_i} W^\T f_{ij}.
\end{equation}
Similar to the matrix-vector multiplication in \cref{sec:matvec} of the nonstandard form, we propose
a neural network for meta-learning \cref{eq:WHWf} with four modules:
\begin{enumerate}
\item $\eta\to S_{\eta}$: a module learns the map $\eta\to C_{\eta}$ and then generates the
  banded sparse matrix $S_{\eta}$ from $C_{\eta}$ (denoted as $S_{\eta} = \mathcal{S}[C_{\eta}]$);
\item $v\sps{L}\to \hat{v}:=W^\T v\sps{L}$: a module applies the forward wavelet transform to
  $v\sps{L}$ to generate $\hat{v}$;
\item $\hat{u} := S_{\eta}\hat{v}$: a module evaluates the matrix-vector multiplication in the
    nonstandard form;
\item $\hat{u}\to u\sps{L}:=W\hat{u}$: a module applies the inverse wavelet transform on
  $\hat{u}$ to generate $u\sps{L}$.
\end{enumerate}

Instead of computing $C_\eta$ from the full operator $G_\eta$ as described in \cref{sec:matvec}.
the first module forms $C_\eta$ directly from the parameter $\eta$ using a deep NN. This module can
be split into two steps: (1) carrying out the map $\eta\to C_{\eta}\sps{\ell}$ for each scale
$\ell$; (2) constructing the nonstandard form $S_{\eta}$ from
$C_{\eta}:=\{C_{\eta}\sps{\ell}\}_{\ell=L_0,\dots,L-1}$. The NN architecture for the map $\eta\to
C_{\eta}\sps{\ell}$ is often problem-dependent. For many applications, including the ones to be
considered in \cref{sec:sch,sec:rte}, the problem is often translation-invariant, i.e., for any
translation operator $T$,
\begin{equation}
  u=G_{\eta}f \quad \text{implies} \quad (T u) = G_{(T\eta)} (Tf).    
\end{equation}
For such problems, a convolutional NN is often used for its efficiency and robustness.

The second and fourth modules perform the forward and inverse wavelet transforms (as in
\cref{sec:matvec}), respectively, for a specific wavelet basis. The selection of an effective
wavelet basis is often problem-dependent. The capability of learning a problem-dependent wavelet
transform from data is essential for the accuracy of the NN architecture.

\begin{algorithm}[htb]
\begin{small}
\begin{center}
\begin{algorithmic}[1]
    \Require $f\sps{L}=f$, $\eta\sps{L}=\eta$, $\alpha_1,\alpha_2\in\bbN$; \Ensure $u$;
    \LineComment{ {\bf Module 1:} Learn the map $\eta\to S_{\eta}$} 
    \For {$\ell$ from $L-1$ to $L_0$ by $-1$}
                \State $C_{\eta}\sps{\ell}=\ConvNet[\ell,\alpha_1,\K](\eta)$
    \EndFor
        \State Generate $D_j\sps{\ell}$, $j=1,2,3$, $\ell=L_0,\dots,L-1$ and $A\sps{L_0}$ from 
    $C_{\eta}\sps{\ell}$, $\ell=L_0,\dots,L-1$ \label{lin:generate}
                                                        \vspace{1mm}
    \LineComment{ {\bf Module 2:} forward wavelet transform on $v\sps{L}$}
    \State $v\sps{L}=f\sps{L}$
    \For {$\ell$ from $L-1$ to $L_0$ by $-1$}
            \State $(d\sps{\ell},v\sps{\ell}) = \FWT[\alpha_2](v\sps{\ell+1})$
        \EndFor
    \vspace{1mm}
    \LineComment{ {\bf Module 3:} matrix-vector multiplication in the nonstandard form}
    \State $D_4\sps{\ell}=0$ for $\ell = L_0+1,\cdots, L-1$ and $D_4\sps{L_0}=A\sps{L_0}$
    \For {$\ell$ from $L_0$ to $L-1$}\vspace{1mm}
    \State 
        $\begin{pmatrix}
            w\sps{\ell}\\
            s\sps{\ell}
        \end{pmatrix}=\begin{pmatrix}
            D_1\sps{\ell} & D_2\sps{\ell}\\
            D_3\sps{\ell} &  D_4\sps{\ell}
        \end{pmatrix}\begin{pmatrix}
            \cw\sps{\ell}\\
            \cs\sps{\ell}
        \end{pmatrix}$
        \label{lin:matvec}
    \EndFor
                                                    \vspace{1mm}
    \LineComment{ {\bf Module 4:} Inverse wavelet transform}
    \State $u\sps{L_0}=0$
    \For {$\ell$ from $L_0$ to $L-1$}
                \State $u\sps{\ell+1}=\IWT[\alpha_2]([w\sps{\ell}, s\sps{\ell}+u\sps{\ell}])$
            \EndFor
                \State Average over the channel direction of $u\sps{L}$ to give $u$;
    \State $\mathsf{return}$ $u$
\end{algorithmic}
\end{center}
\end{small}
\caption{Neural network architecture for meta-learning  $u=G_{\eta} f$.}
\label{alg:NN}\end{algorithm}

\begin{figure}[h!]
  \centering
  \includegraphics[width=1\textwidth,trim=2cm 1cm 1.5cm 0cm, clip]{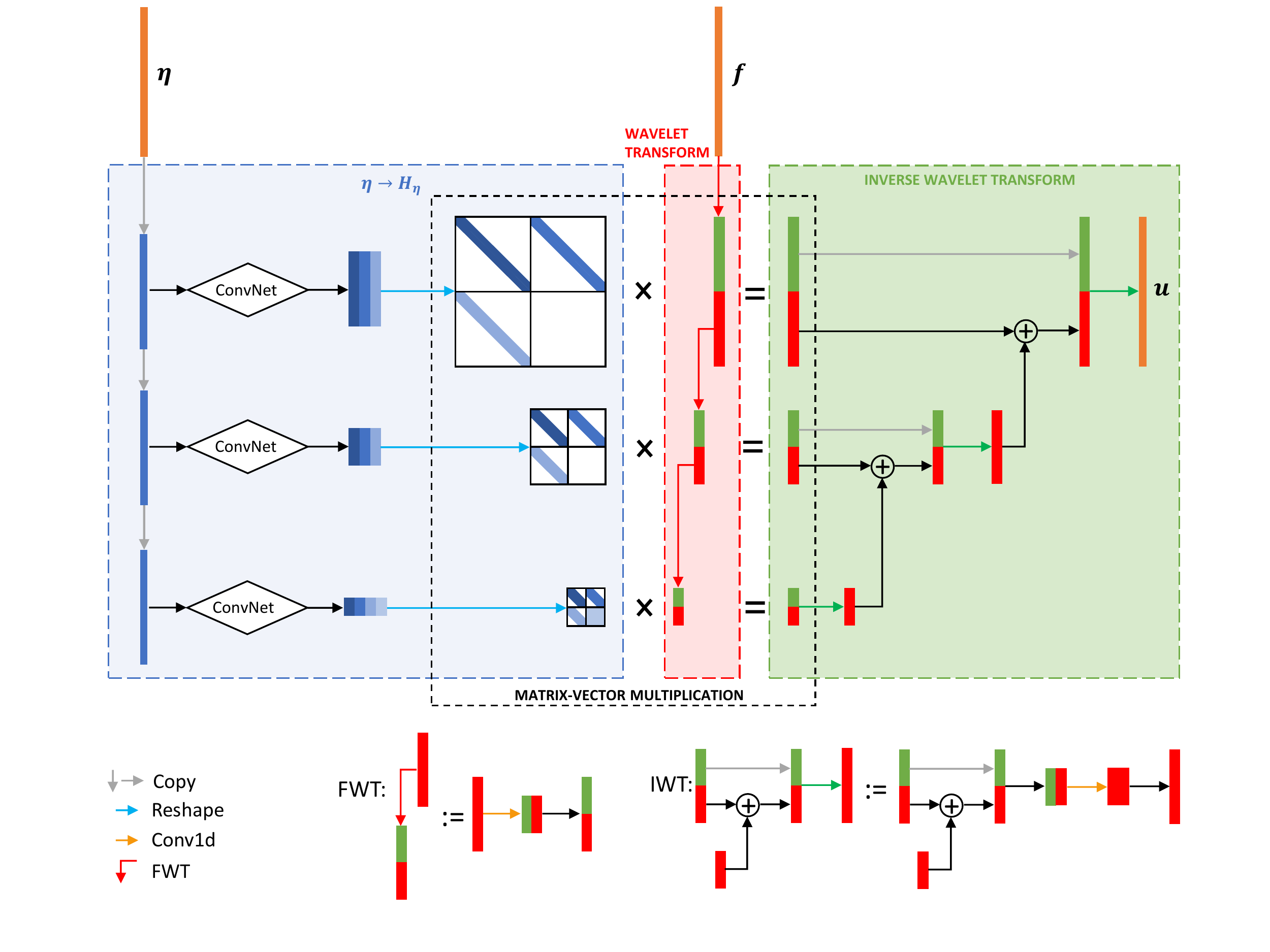}
  \caption{\label{fig:NN_arch} Illustration of the neural network architecture from
    \cref{alg:NN}, with input vectors $\eta$ and $f$ and output vector $u$. Each of the four
      modules in \cref{alg:NN} is represented by blocks with dashed contours. From left to right,
      the map $\eta\to S_{\eta}$ is the blue block, the forward wavelet transform applied to $f$
      is the red block, the sparse matrix-vector multiplication with $S_\eta$ is the transparent
      block, and the inverse wavelet transform is the green block.
      }
\end{figure}

Combining these four modules results in the architecture summarized in \cref{alg:NN}. An
illustration is given in \cref{fig:NN_arch}. Below we describe details of the layers and parameters
used in this architecture.

\begin{figure}[h!]
  \centering
  \includegraphics[width=0.7\textwidth,page=2]{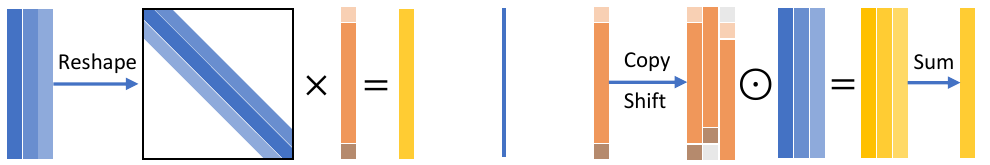}
  \caption{\label{fig:diagonal_matvec}Two approaches to evaluate the multiplication of a
    band matrix (stored with vectors) with a vector. 
    The left figure corresponds to the case in \cref{alg:NN} and the right one is used in the 
    implementation. The symbol $\odot$ stands for element-wise multiplication.
  }
\end{figure}

\paragraph{Implementation details.} 
The input, output, and intermediate data of the NN architecture are all represented with
$2$-tensors. For a tensor of size $N\times \alpha$, we refer to $N$ as the spatial dimensions and
$\alpha$ as the channel dimension.  The main tool is the \emph{convolutional layer}. Given an input
tensor $\xi$ of size $N\times\alpha$, the convolutional layers outputs a tensor $\zeta$ of size
$N'\times \alpha'$ obtained via
\begin{equation}\label{eq:conv}
    \zeta_{i,c'} = \phi\left( \sum_{j=is}^{is+w-1}\sum_{c=0}^{\alpha-1}
    W_{j;c',c}\xi_{j,c} + b_{c'}\right),
    \quad i=0,\dots,N'-1, ~ c'=0,\dots,\alpha'-1,
\end{equation}
where $w$ is the \emph{window size}, $s$ is the \emph{stride} and $\phi$ is the \emph{activation
function}, usually chosen to be a linear function, a rectified-linear unit (ReLU) function, or a
sigmoid function. We denote this convolutional layer as
\begin{equation}\label{eq:conv_compact}
    \zeta = \Convone[\alpha', w, s, \phi](\xi).
\end{equation}
The basic building blocks and layers used in \cref{alg:NN} are listed below.
\begin{itemize}

\item $\eta \to C_{\eta}\sps{\ell}:C_{\eta}\sps{\ell}=\ConvNet[\ell,\alpha,\K](\eta)$. As discussed
  above, it is often a convolutional NN if the system \cref{eq:Gf} is translation invariant. Since
  the spatial size of $\eta$ is greater than that of $C_{\eta}\sps{\ell}$,
  $\ConvNet[\ell,\alpha,\K]$ consists of $\K$ convolutional layers and several downsampling or
  pooling layers.

                  \item Forward wavelet transform at level $\ell$: $(d\sps{\ell}, v\sps{\ell}) =
  \FWT[\alpha](v\sps{\ell+1})$. This is the NN representation of the first equation in
  \cref{eq:FWT}. It is implemented as $f\sps{\ell} = \Convone[2\alpha, 2p, 2, \id](v\sps{\ell})$,
  where the first $\alpha$ and the last $\alpha$ channels of $f\sps{\ell}$ are assigned to
  $d\sps{\ell}$ and $v\sps{\ell}$, respectively.
\item Inverse wavelet transform at level $\ell$: $u\sps{\ell+1} = \IWT[\alpha]([w\sps{\ell},
  s\sps{\ell}])$. This is the NN representation of the second equation in \cref{eq:FWT}. The
  expression $[w\sps{\ell}, s\sps{\ell}]$ stands for concatenating the $2$-tensors $w\sps{\ell}$ and
  $s\sps{\ell}$ of size $2\sps{\ell}\times \alpha$ to a $2$-tensor of size $2\sps{\ell}\times
  2\alpha$ along the channel dimension. This layer first applies the inverse transform, implemented
  by $\Convone[2\alpha, p, 1, \id]$, and then reshapes the output of size $2^{\ell}\times 2\alpha$
  to a $2$-tensor of size $2^{\ell+1}\times \alpha$ by a column-first ordering.
\end{itemize}
The generation of $D_{j}\sps{\ell}$ and $A\sps{L_0}$ from $C_{\eta}^{\ell}$ in \cref{lin:generate}
and the matrix-vector multiplication in \cref{lin:matvec} of \cref{alg:NN} require some
discussion. \Cref{fig:diagonal_matvec} illustrates two approaches for evaluating the matrix-vector
multiplication of a band matrix whose nonzero entries are stored in a set of vectors. The left
figure corresponds to the case in \cref{alg:NN}, while the right one is used in the actual
implementation. To avoid the copying and shifting of $d\sps{\ell}$ and $v\sps{\ell}$, it is
convenient to set $\alpha_2=\alpha_1=\alpha$. Though there are slightly more NN parameters in this
case, this implementation change allows for a more flexible NN that can learn faster.

\subsection{The multidimensional case}
Let us focus on the 2D case. The input, the output and the intermediate data are all $3$-tensors of
size $N_1\times N_2\times \alpha$, where $N = (N_1,N_2)$ is the spatial dimension and $\alpha$ is
the channel dimension. The convolutional layer takes the form
\begin{equation}\label{eq:conv2}
  \zeta_{i,c'} = \phi\left( \sum_{j_1=i_1s}^{i_1s+w-1}\sum_{j_2=i_2s}^{i_2s+w-1}\sum_{c=0}^{\alpha-1}
  W_{j;c',c}\xi_{j,c} + b_{c'}\right),
  \quad i_1=0,\dots,N_1'-1, i_2=0, \dots,N_2'-1, c'=0,\dots,\alpha'-1,
\end{equation}
where the input $\xi$ is of size $N_1\times N_2\times\alpha$ and the output $\zeta$ is of size
$N_1'\times N_2'\times\alpha'$.  Here, the same stride $s$ and window size $w$ are used in both
dimensions. We denote this convolutional layer as
\begin{equation}
  \zeta= \Convtwo[\alpha, w, s, \phi](\xi).
\end{equation}
\cref{alg:NN} can be easily extended to the 2D case, following the same way that \cref{alg:BCR} was
extended in \cref{sec:matvec_2d}. Since there are three different types of wavelets in 2D, the
layers in \cref{alg:NN} are redefined as follows:
\begin{itemize}
\item   $\eta\to C_{\eta}\sps{\ell}$ module: $C_{\eta}\sps{\ell}=\ConvNet[\ell,\alpha,\K](\eta)$. This
  module is often a two-dimensional convolutional NN with several downsampling or pooling layers.
\item Wavelet transform at level $\ell$: $(d_1\sps{\ell}, d_2\sps{\ell}, d_3\sps{\ell}, v\sps{\ell})
  = \FWT[\alpha](v\sps{\ell+1})$. This is implemented using $f\sps{\ell} = \Convtwo[4\alpha, 2p, 2,
    \id](v\sps{\ell})$. The first, second, third and last $\alpha$ channels of $f\sps{\ell}$ are
  assigned to $d_1\sps{\ell}$, $d_2\sps{\ell}$, $d_3\sps{\ell}$ and $v\sps{\ell}$, respectively.
\item Inverse wavelet transform at level $\ell$: $u\sps{\ell+1} = \IWT[\alpha]([w_1\sps{\ell},
  w_2\sps{\ell}, w_3\sps{\ell}, s\sps{\ell}]$. This is implemented by first computing
  $\Convtwo[4\alpha, p, 1, \id]([d_1\sps{\ell}, d_2\sps{\ell}, d_3\sps{\ell},
  v\sps{\ell}+u\sps{\ell}])$, and then reshaping the output of size $2^{\ell}\times 2^{\ell}\times
  4\alpha$ to a $3$-tensor of size $2^{\ell+1}\times 2^{\ell+1}\times \alpha$. The reshape operation
  is performed as follows: (1) reshape the output to a $5$-tensor of size
  $2^\ell\times 2^\ell\times 2\times 2\times \alpha$ by splitting the last dimension; (2) permute
  the second and third dimensions to obtain a $5$-tensor of size $2^\ell\times 2\times 2^\ell\times
  2\times \alpha$; (3) group the first and second dimensions, and the third and fourth dimensions,
  respectively, to obtain the resulting $3$-tensor of size $2^{\ell+1}\times 2^{\ell+1}\times
  \alpha$.
\end{itemize}

\section{Elliptic partial differential equations}\label{sec:sch}
This section applies the meta-learning approach described in \cref{sec:nn} to the Green's functions
of elliptic PDEs, both in the Schr{\"o}dinger form and in the divergence form.

\subsection{Schr{\"o}dinger form}

Consider the equation
\begin{equation}\label{eq:nlse}
  \begin{aligned}
    &-\Delta u(x) + \eta(x)u(x) =f(x),\quad x\in \Omega = [0,1]^d,\\
  \end{aligned}
\end{equation}
with a periodic boundary condition, where $\eta(x)>0$ is the potential and $f(x)$ is the source
term. Following the notations of \cref{sec:intro},
\begin{equation}\label{eq:schrodinger_Gf}
  L_\eta = -\Delta + \eta(x), \qquad G_{\eta}=L_\eta^{-1}, \qquad     u = G_{\eta}f.
\end{equation}
Since the problem \cref{eq:nlse} is translation-invariant due to the periodic boundary condition,
the map $\eta\to C_{\eta}\sps{\ell}$ can be represented with a convolutional NN. In what follows, we
first derive the explicit dependence of $C_{\eta}\sps{\ell}$ on $\eta$ using a linear perturbative
analysis and then report some numerical studies.

\paragraph{Mathematical analysis.} When $\eta$ is close to a fixed homogeneous background $\eta_0>0$, it is convenient to write
\begin{equation}\label{eq:se_op}
  L_\eta =  L_0-E_{\eta}, \quad L_0=-\Delta+\eta_0,\quad E_{\eta}=\diag(-\eta+\eta_0), 
\end{equation}
Let $G_0=L_0^{-1}$ be the Green's function of $L_0$ with the periodic boundary condition.  Using the
Neumann series for the resolvent $(I-G_0E_{\eta})^{-1}$ with $|\eta(x)-\eta_0|$ sufficiently small,
one can write the Green's function $G_{\eta}$ as a perturbative expansion
\begin{equation}  \label{eq:se_g}
  G_{\eta} = (L_0-E_{\eta})^{-1} = G_0 + G_0 E_{\eta} G_0 + G_0 E_{\eta} G_0 E_{\eta} G_0 +
  \ldots.
\end{equation}
For sufficiently small $|\eta(x)-\eta_0|$, the operator $G_{\eta}$ can be approximated by its linear
part as
\begin{equation}  \label{eq:se_gapprox}
  G_{\eta}\approx G_0 +  G_0 E_{\eta} G_0.
\end{equation}

Let $g_0(x)$ and $g_{\eta}(x)$ be the kernel of $G_0$ and $G_{\eta}$, respectively.
Since $G_0$ is the Green's function of $-\Delta+\eta_0$ with the periodic boundary condition, the
kernel $g_0$ is translation-invariant, i.e., $g_0(x,y)=g_0(x-y)$.
The wavelet-wavelet coefficients of $G_{\eta}$ at level $\ell$ take the form 
\begin{equation} \label{eq:se_d1}
  \begin{aligned} 
    D\wl{\ell}{1,k_1,k_2}&=\iint \psi\wl{\ell}{k_1}(x)g_{\eta}(x,y)\psi\wl{\ell}{k_2}(y)\dd x\dd y \\
    & \approx 
     \iint \psi\wl{\ell}{k_1}(x)g_0(x-y)\psi\wl{\ell}{k_2}(y)\dd x\dd y
    + \iiint \psi\wl{\ell}{k_1}(x)g_0(x-z)(\eta(z)-\eta_0)g_0(z-y)\psi\wl{\ell}{k_2}(y)\dd x\dd y \dd z\\
    &= \iint \psi\wl{\ell}{k_1}(x)g_0(x-y)\psi\wl{\ell}{k_2}(y)\dd x\dd y
     + \int \widetilde{\psi}\wl{\ell}{k_1}(z) \widetilde{\psi}\wl{\ell}{k_2}(z) (\eta(z)-\eta_0) \dd z,
  \end{aligned}
\end{equation}
where $\widetilde{\psi}\wl{\ell}{k}(z) := (g_0 * \psi\wl{\ell}{k})(z)$. For a fixed diagonal of
$D_1^\ell$ with $k_2=k_1+c$ for a constant $c$, \cref{eq:se_d1} states that the map from $\eta$ to
$D^\ell_{1,k_1,k_1+c}$ for all possible $k_1$ is simply a convolution with an addition of a term
independent of $\eta$, which can be simply represented by the \Convone layer in
\cref{eq:conv}.
It is straightforward to extend the conclusion to $D\sps{\ell}_{j,k_1,k_2}$, $j=2,3$ and
$A\sps{L}_{k_1,k_2}$.

When $|\eta(x)-\eta_0|$ is not small, one can account for the nonlinearities neglected in the
perturbative analysis by using multiple convolutional layers and making use of nonlinear activation
functions. In other words, it is natural to approximate the map $\eta\to C_{\eta}\sps{\ell}$ using a
convolutional NN with enough layers and an appropriate window size \cite{Liang2017,Jacot2018,Ilsang2019}.

Moreover, since the matrix $G_{\eta}$ is a symmetric matrix, $D_1\sps{\ell}$ and $A\sps{L_0}$ are
symmetric and $(D_2\sps{\ell})^\T=D_3\sps{\ell}$. In the implementation, the symmetry is enforced by
generating $D_3\sps{\ell}$ from $D_2\sps{\ell}$, and replacing $D_1\sps{\ell}$ (or $A\sps{L_0}$) by
$\frac{1}{2}(D_1\sps{\ell}+(D_1\sps{\ell})^\T)$ (or $\frac{1}{2} (A\sps{L_0} + (A\sps{L_0})^\T)$,
respectively.  Since the Schr\"odinger form considered in this section includes the periodic
boundary condition, the convolutional layers are all implemented with periodic padding.

\paragraph{Numerical results.}
The NN discussed above is implemented in Keras \cite{keras} (running on top of TensorFlow
\cite{tensorflow}). The parameters of the NN are initialized randomly from the normal
distribution. The loss function is set to be the mean squared error
\begin{equation}
    \frac{1}{\Nsamp}\|u-u^{\NN}\|_{\ell^2},
\end{equation}
where the exact solution, obtained by solving \cref{eq:schrodinger_Gf}, is denoted as $u$ and the NN
prediction as $u^{\NN}$. $\Nsamp$ denotes the number of samples. The NN is trained until convergence
using the Nadam optimizer \cite{dozat2016incorporating} with the learning rates equal to $10^{-3}$
for the 1D case and $10^{-4}$ for the 2D case. The batch size is set to be one percent of the number
of training samples. The support of the scaling function $\varphi$ is chosen to be $2p=6$. The
number of levels $L-L_0$ in the wavelet transform is $6$ for the 1D case and $4$ for the 2D case.

The data set contains $5,000$ different $\eta$ and for each $\eta$ \cref{eq:schrodinger_Gf} is
solved with $20$ randomly generated $f$ using the central difference scheme. Therefore, the number of training
samples corresponds to the number of different $\{f_{ij}\}$, rather than different
$\{\eta_i\}$. Half of the generated data is used for training data, while the other half is reserved
for testing. The accuracy of the NN is measured by the relative error in
the $\ell^2$ norm
\begin{equation}
  \epsilon = \frac{\|u-u^{\NN}\|_{\ell^2}}{\|u\|_{\ell^2}}.
\end{equation}
The training error $\trainerror$ and test error $\testerror$ are calculated by averaging the
relative error over all training and test samples, respectively. The number of parameters in the NN
is denoted by $\Nparams$. The operator error $\operror$ is calculated by averaging the
  relative $2$-norm error of the matrix
\begin{equation}
  \frac{\|G_\eta-G_\eta^{\NN}\|_{\ell^2\rightarrow\ell^2}}{\|G_\eta\|_{\ell^2\rightarrow\ell^2}}
\end{equation}
over samples of the exact inverse operator $G_\eta$ and its NN approximation $G_\eta^{\NN}$.

\begin{table}[h!]
  \centering
  \begin{tabular}{cccccc}
    \hline
    $\alpha$ & $\K$ & $\Nparams$ & $\trainerror$ & $\testerror$ & \RRev{$\operror$} \\ \hline\hline
    5 & 5  &  30201  &  \RRev{4.43e-3}   &    \RRev{4.74e-3}  & \RRev{2.49e-3} 	\\ \hline
    5 & 7  &  38061 & \RRev{4.83e-3}  &  \RRev{5.18e-3} & \RRev{4.28e-3}  \\ \hline
    7 & 5 &   58717 & \RRev{4.09e-3}    &   \RRev{4.35e-3} & \RRev{3.28e-3} \\ \hline
    7 & 7 &   74089  &  \RRev{4.11e-3}  &  \RRev{4.42e-3} & \RRev{2.18e-3}  \\ \hline
  \end{tabular}
  \caption{\label{tab:se1d} Relative error in approximating the solution to the Schr\"odinger form in the 1D case.}
\end{table}

\begin{figure}[h!]
  \centering
  \subfloat[\label{fig:se1d_input} $\eta$]{
    \includegraphics[width=0.4\textwidth,trim=0.55cm 1.15cm 0.75cm 0.95cm,  clip]{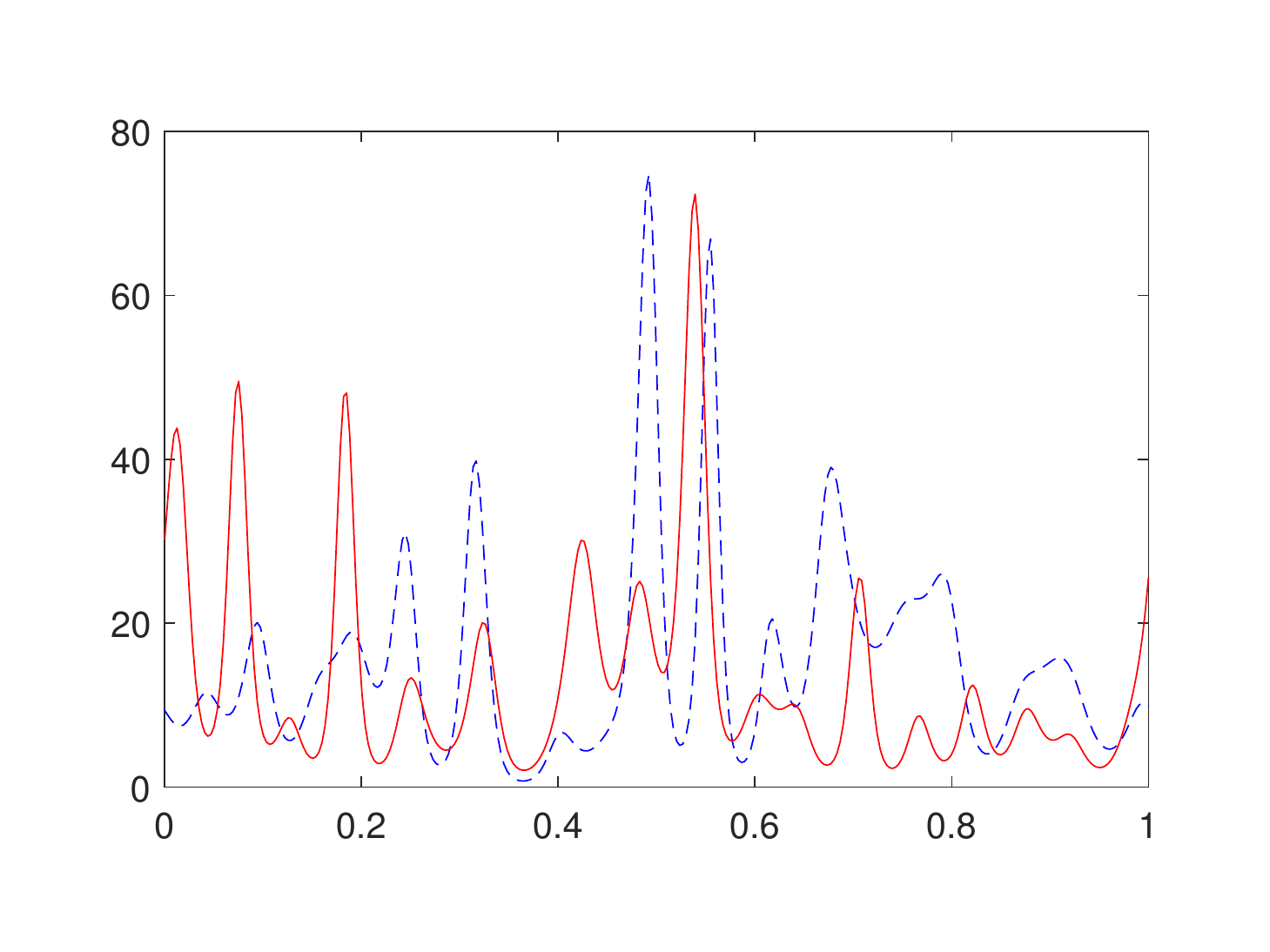}
  }
  \subfloat[\label{fig:se1d_input_f} $f$]{
    \includegraphics[width=0.4\textwidth,trim=0.55cm 1.15cm 0.75cm 0.95cm,  clip]{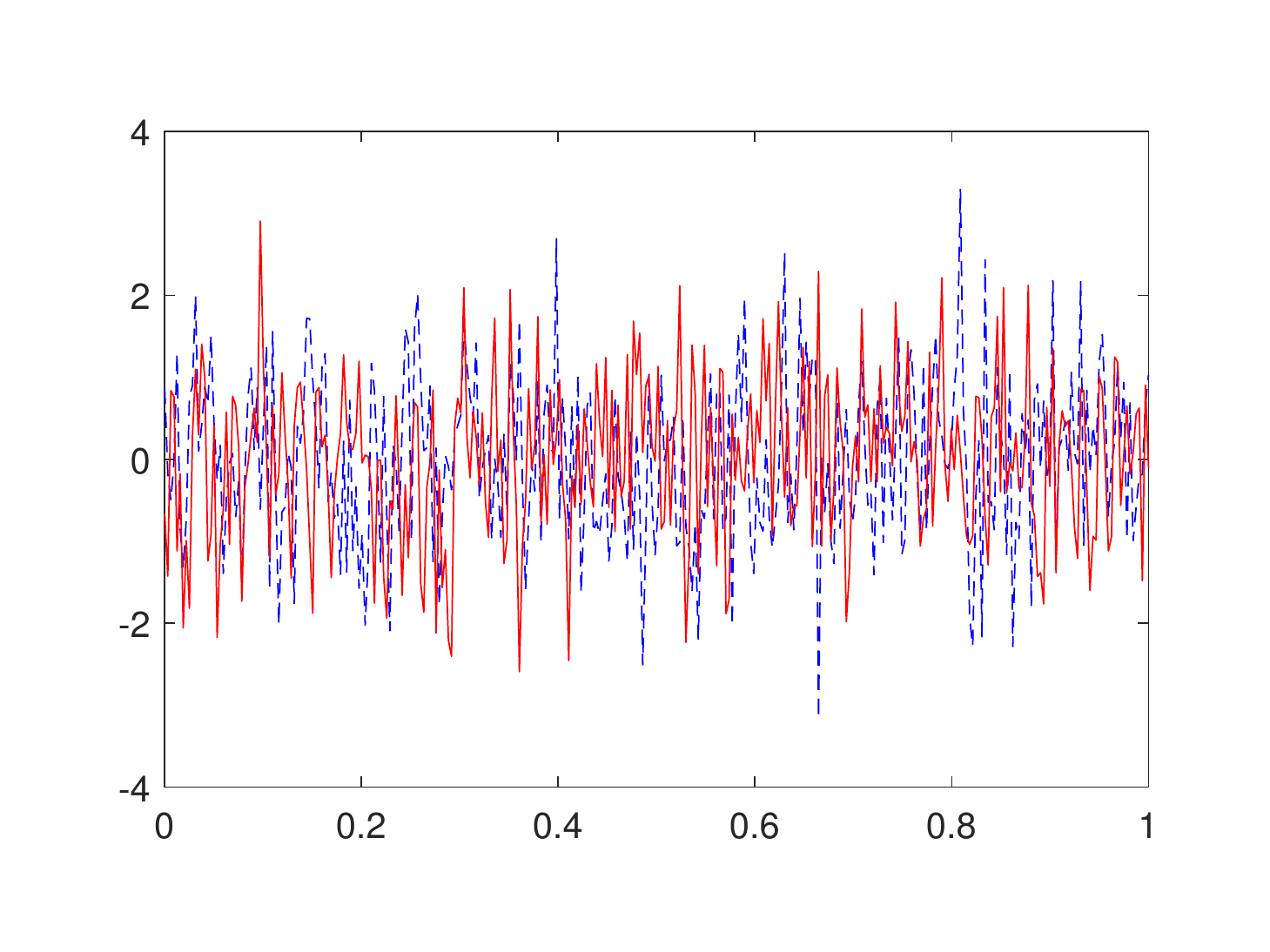}
  }\\
  \subfloat[\label{fig:se1d__out} $u^{\NN}$]{
    \includegraphics[width=0.4\textwidth,trim=0.55cm 1.15cm 0.75cm 0.95cm,  clip]{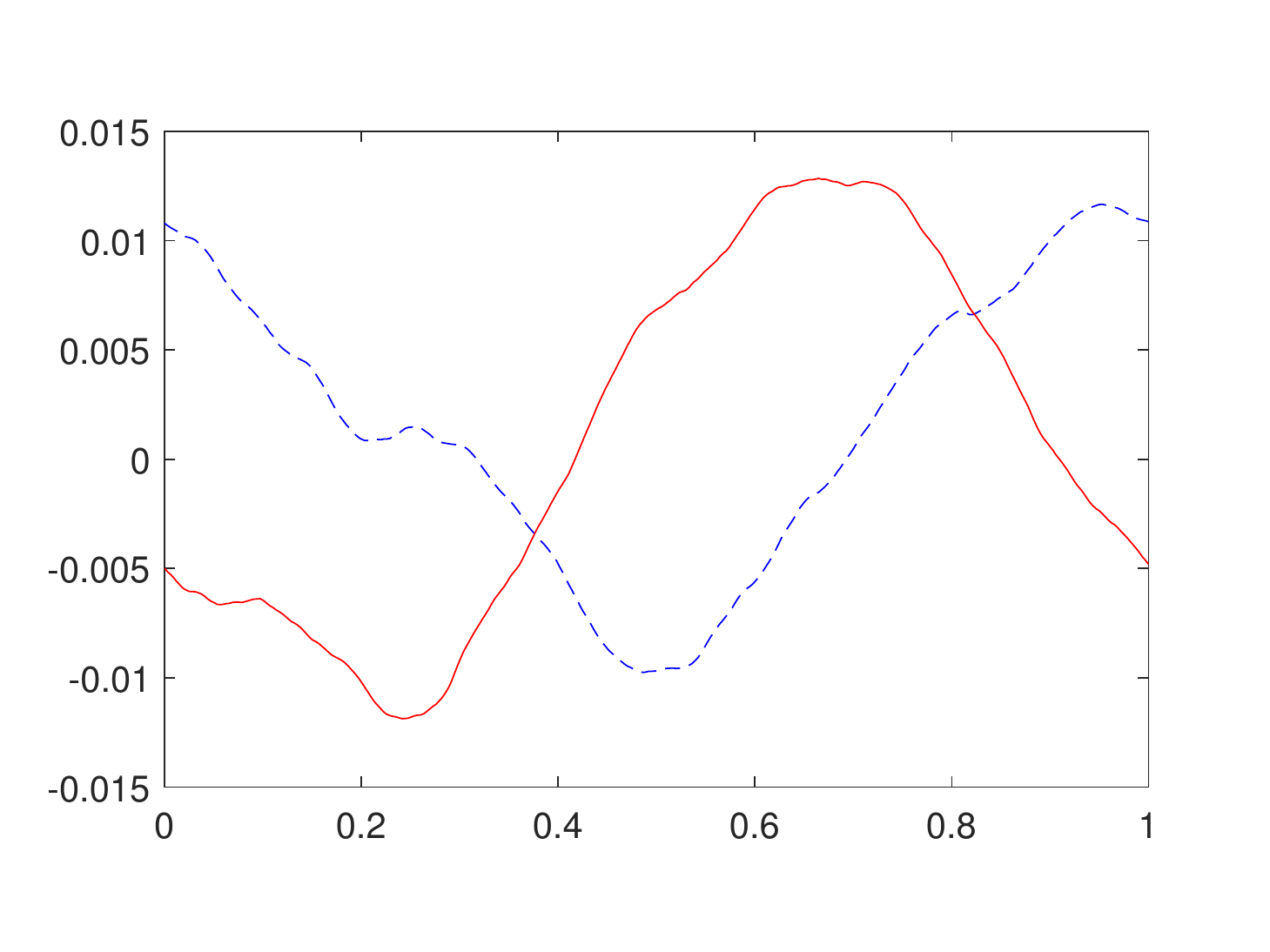}
  }
  \subfloat[\label{fig:se1d_err} $u^{\NN}-u$]{
    \includegraphics[width=0.4\textwidth,trim=0.55cm 1.15cm 0.75cm 0.95cm,  clip]{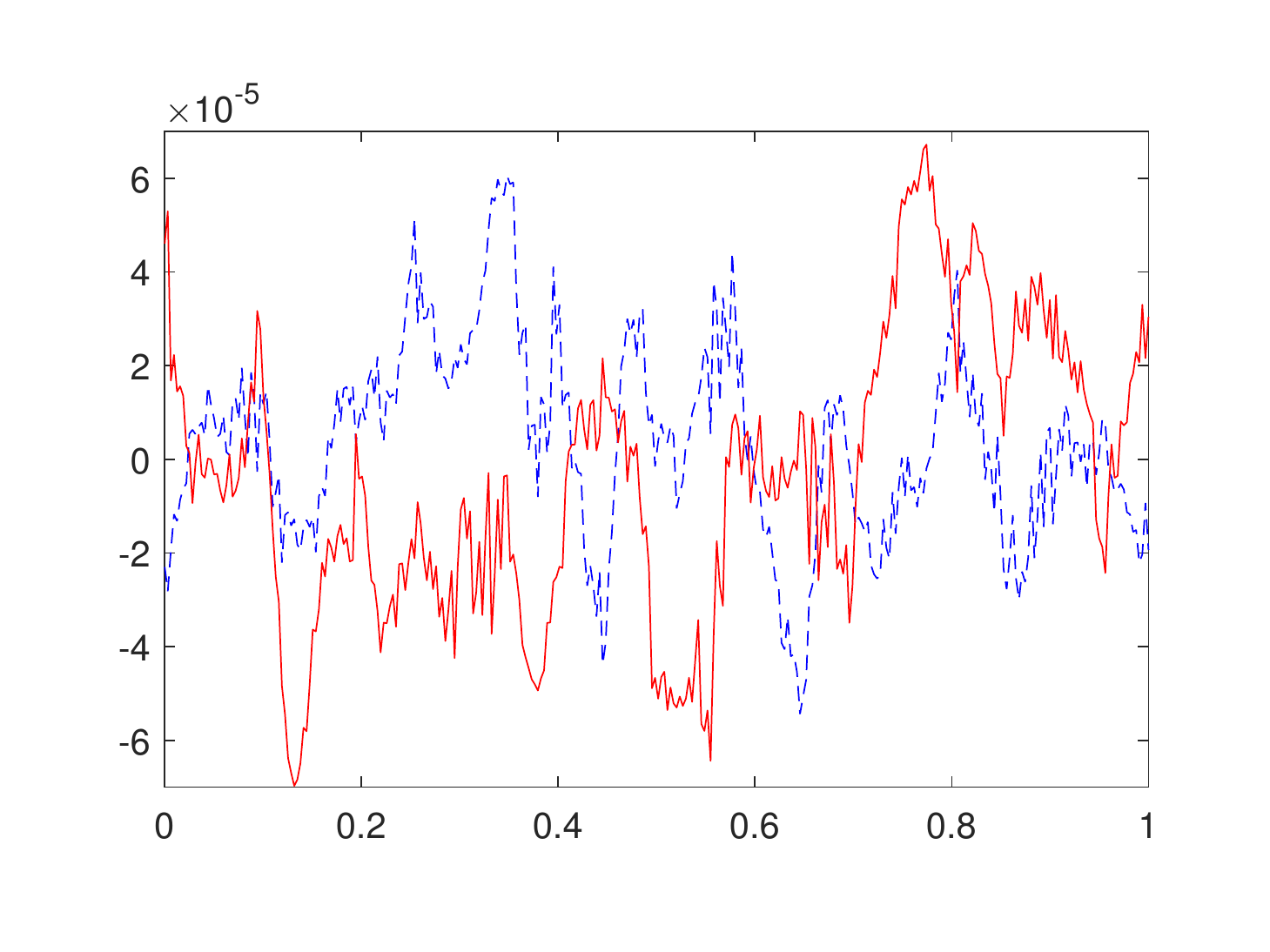}
  }
  \caption{\label{fig:se1d} Two samples (one in red line and the other in blue dashed line) from the
  test set for the potentials $\eta$, the source terms $f$, the predictions $u^{\NN}$ with $\alpha=5$
  and $\K=5$, and their corresponding error for the Schr\"odinger form in the 1D case.}
\end{figure}

\begin{figure}[h!]
  \centering
  \subfloat[\label{fig:se1d_op_ref} $G_\eta$]{
    \includegraphics[width=0.3\textwidth,trim=2.05cm 0.95cm 1.1cm 0.1cm,  clip]{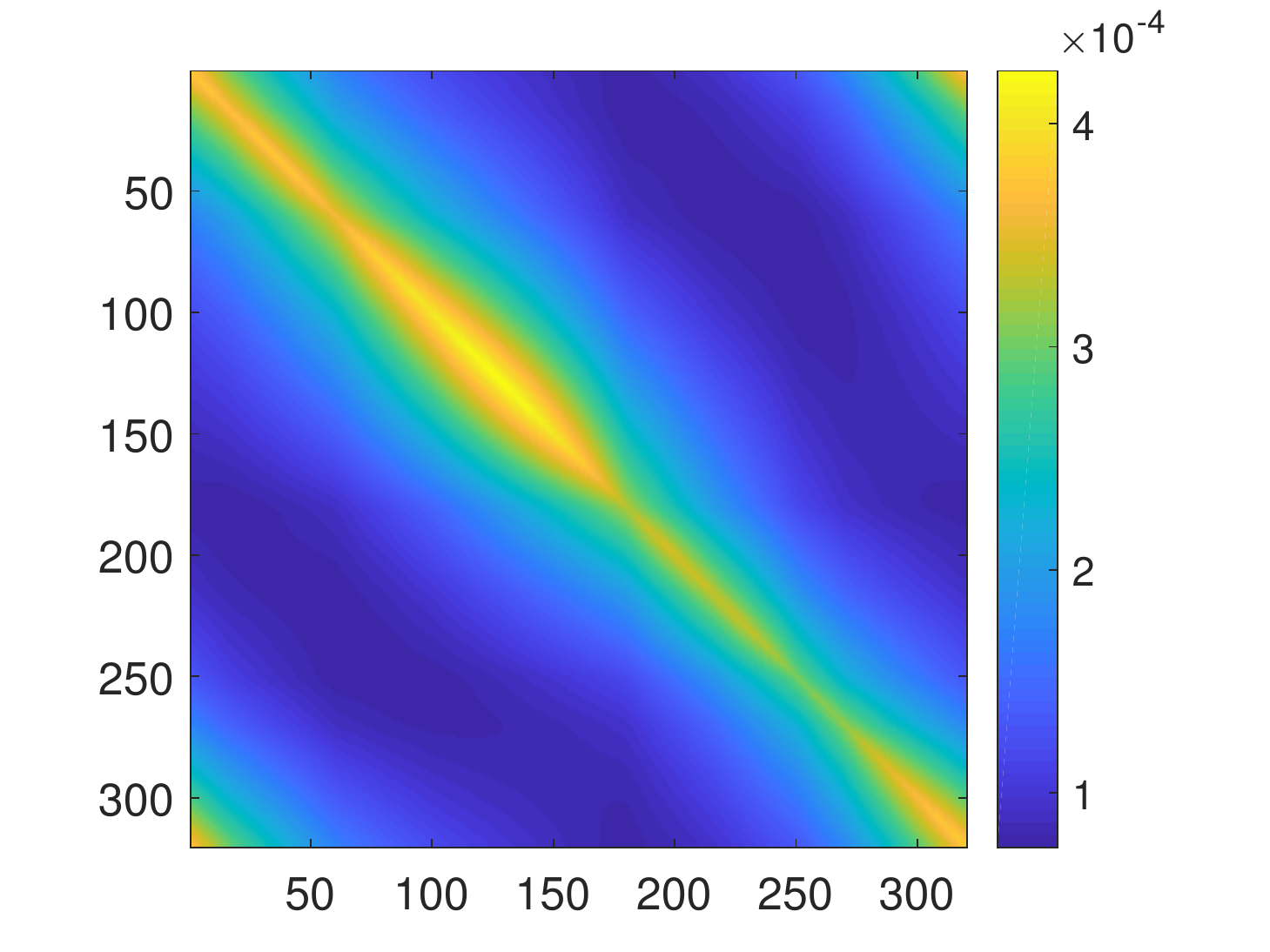}
  }
  \subfloat[\label{fig:se1d_op2} $G_\eta^{\NN}$]{
    \includegraphics[width=0.3\textwidth,trim=2.05cm 0.95cm 1.1cm 0.1cm,  clip]{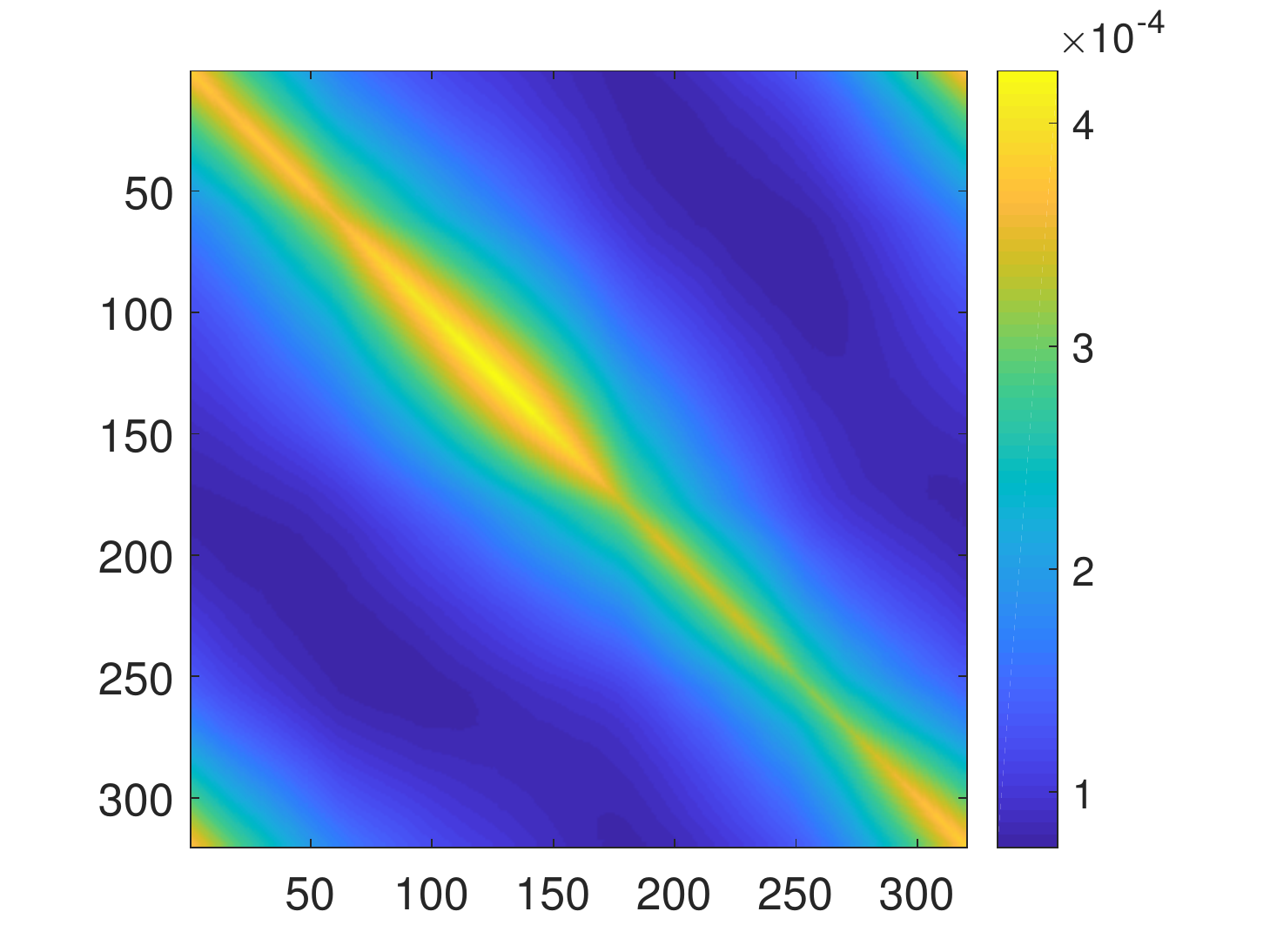}
  }
  \subfloat[\label{fig:se1d_op_dif} $G_\eta^{\NN}-G_\eta$]{
    \includegraphics[width=0.3\textwidth,trim=2.05cm 0.95cm 1.1cm 0.1cm,  clip]{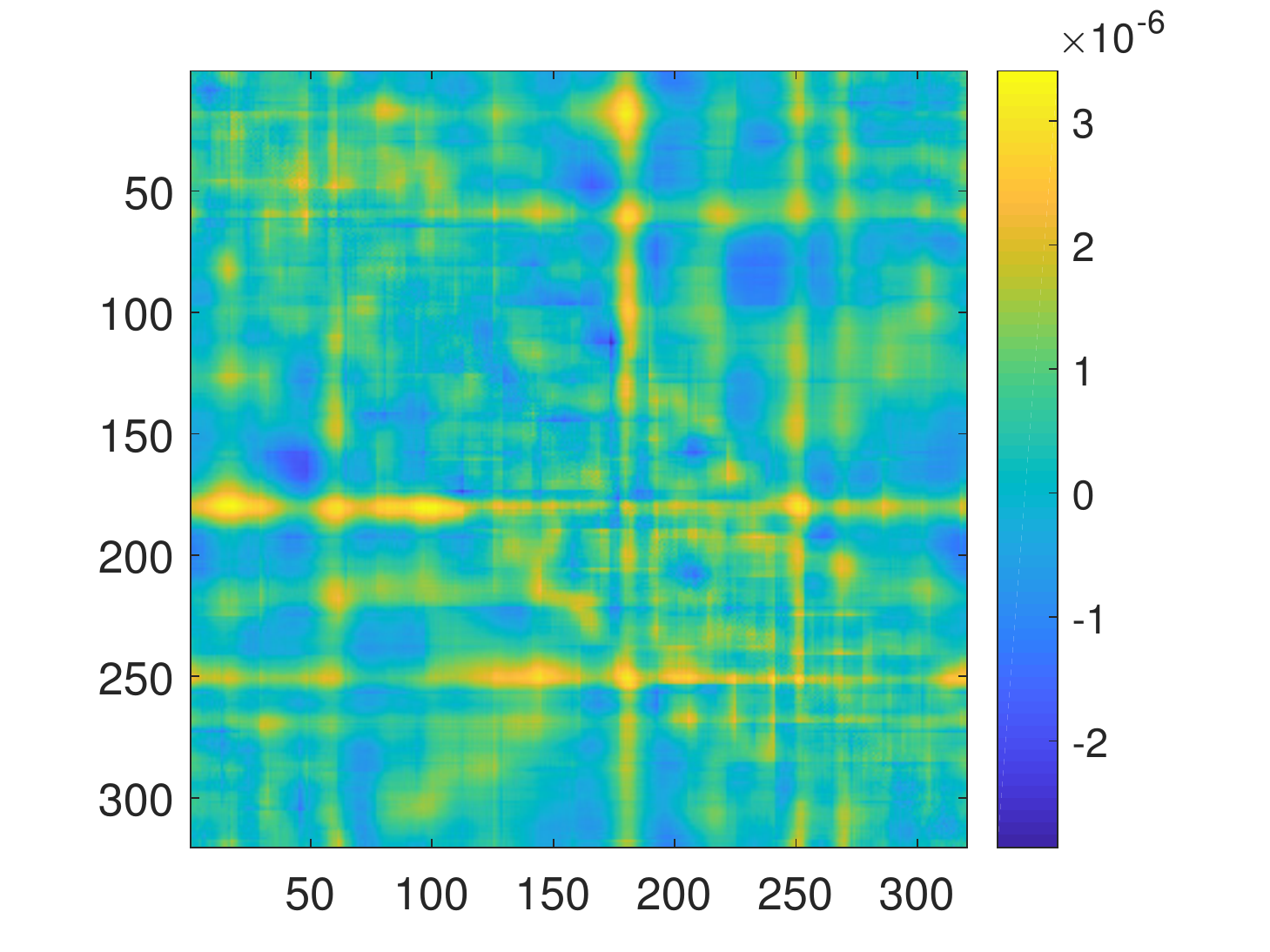}
  }
  \caption{\label{fig:se1d_op} \RRev{Operator approximation with $\alpha=5$
  and $\K=5$ for the Schr\"odinger form in the 1D case.}}
\end{figure}

For the 1D case, the domain $\Omega=[0,1]$ is discretized by a uniform Cartesian grid with $320$
points. The positive potential $\eta(x)$ is generated by (1) sampling independently from
$\mathcal{N}(0,1)$ on a uniform grid with $40$ points, (2) interpolating to the $320$-point grid via
a Fourier interpolation, and (3) point-wise exponentiating followed by a factor of 10
scaling. The source term $f(x)$ is generated by sampling independently from $\mathcal{N}(0,1)$.
The results for different values of $\alpha$ (channel number) and $\K$
(layer number) are reported in \cref{tab:se1d}. The best approximation \RRev{of the operator}, obtained with $\alpha=5$ and
$\K=5$, results in a test error of \RRev{$4.7\times10^{-3}$ and an operator error of $2.5\times10^{-3}$} with only $3\times 10^4$ parameters. \RRev{The operator error reported in \cref{tab:se1d} has been averaged among 100 different samples of $G_\eta$}. Two
random samples from the test data are illustrated in \cref{fig:se1d} along with the NN prediction. \RRev{A representative sample of the inverse operator $G_\eta$ and its NN approximation are displayed in \cref{fig:se1d_op}} 

\begin{table}[h!]
    \centering
    \begin{tabular}{cccccc}
      \hline
      $\alpha$ & $\K$ & $\Nparams$ & $\trainerror$ & $\testerror$ & \RRev{$\operror$} \\ \hline\hline
      11 & 5 & 930447 & \RRev{ 2.21e-2} & \RRev{2.18e-2} & \RRev{4.17e-3}	\\
      15 & 5 & 1226071 & \RRev{ 2.12e-2} & \RRev{2.10e-2}  & \RRev{2.04e-3}  \\  
      \hline
    \end{tabular}
    \caption{\label{tab:se2d} Relative error in approximating the solution to the Schr\"odinger form in the 2D case.}
\end{table}

\begin{figure}[h!]
  \centering
  \subfloat[\label{fig:se2d_input} $\eta$]{
    \includegraphics[width=0.4\textwidth,trim=0.65cm 0.55cm 0.75cm 0.85cm,  clip]{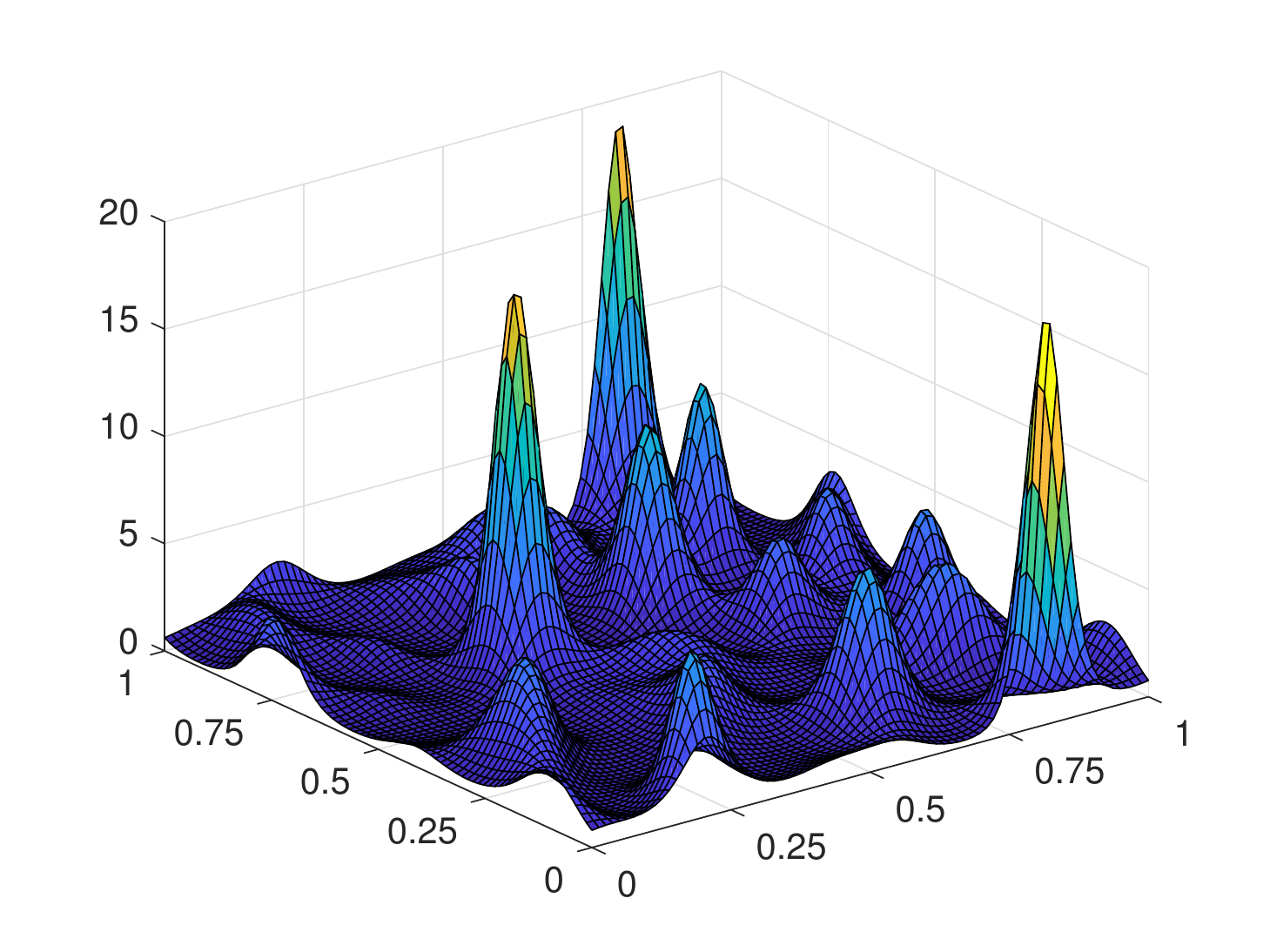}
  }
  \subfloat[\label{fig:se2d_input_f} $f$]{
    \includegraphics[width=0.4\textwidth,trim=0.65cm 0.55cm 0.75cm 0.85cm,  clip]{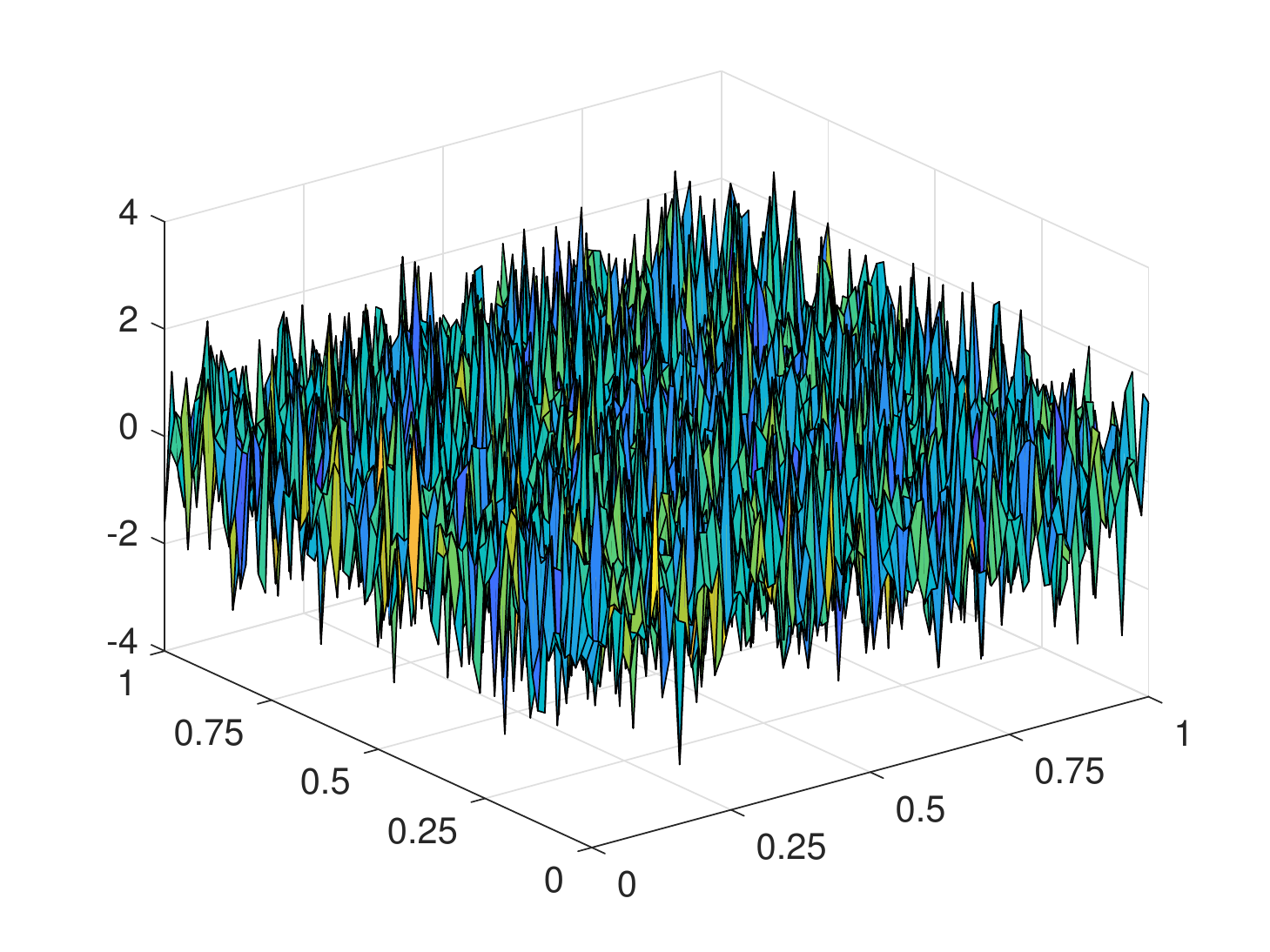}
  }\\
  \subfloat[\label{fig:se2d__out} $u^{\NN}$]{
    \includegraphics[width=0.4\textwidth,trim=0.55cm 0.55cm 0.75cm 0.85cm,  clip]{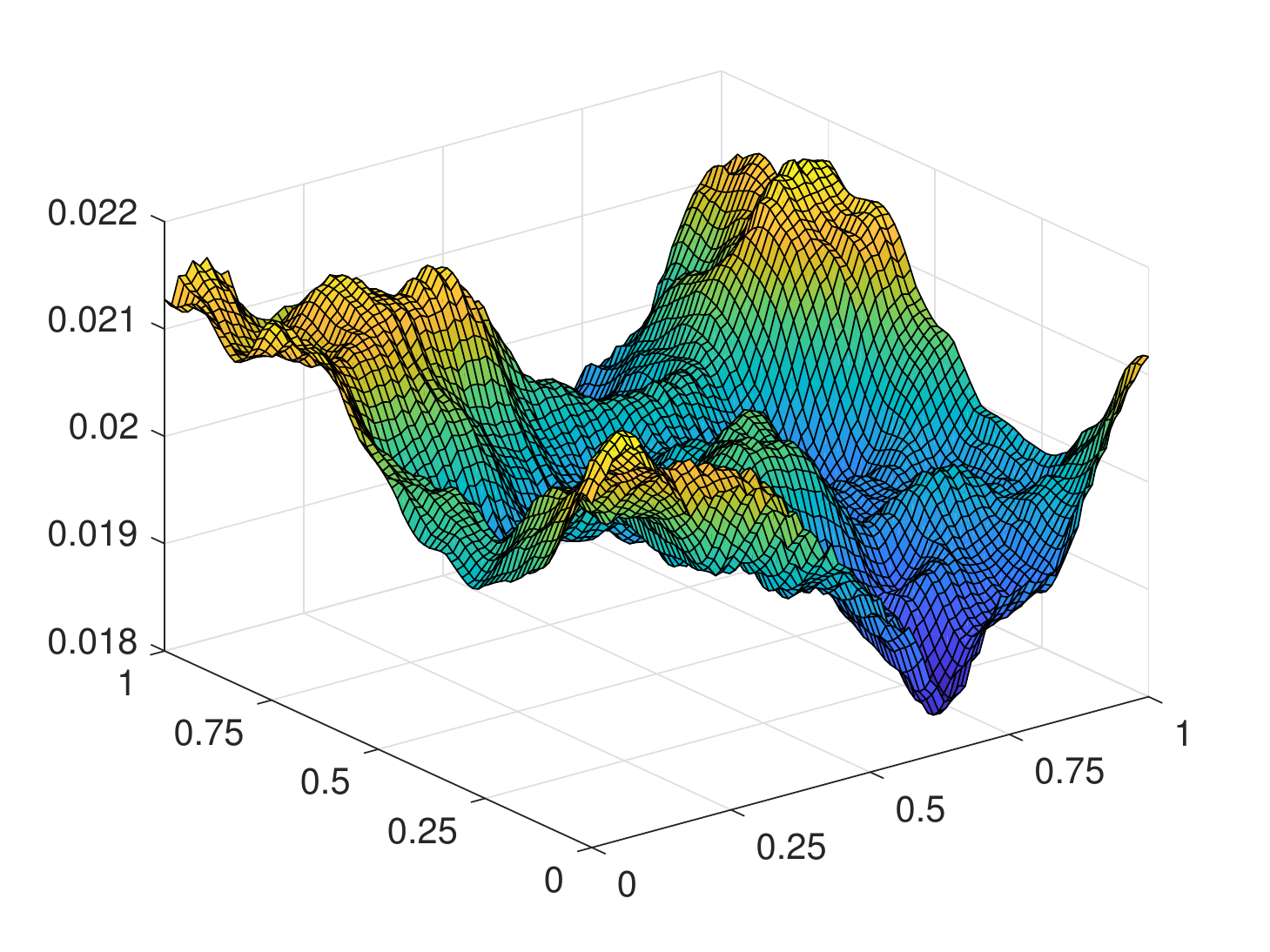}
  }
  \subfloat[\label{fig:se2d_err} $u^{\NN}-u$]{
    \includegraphics[width=0.4\textwidth,trim=0.65cm 0.55cm 0.75cm 0.85cm,  clip]{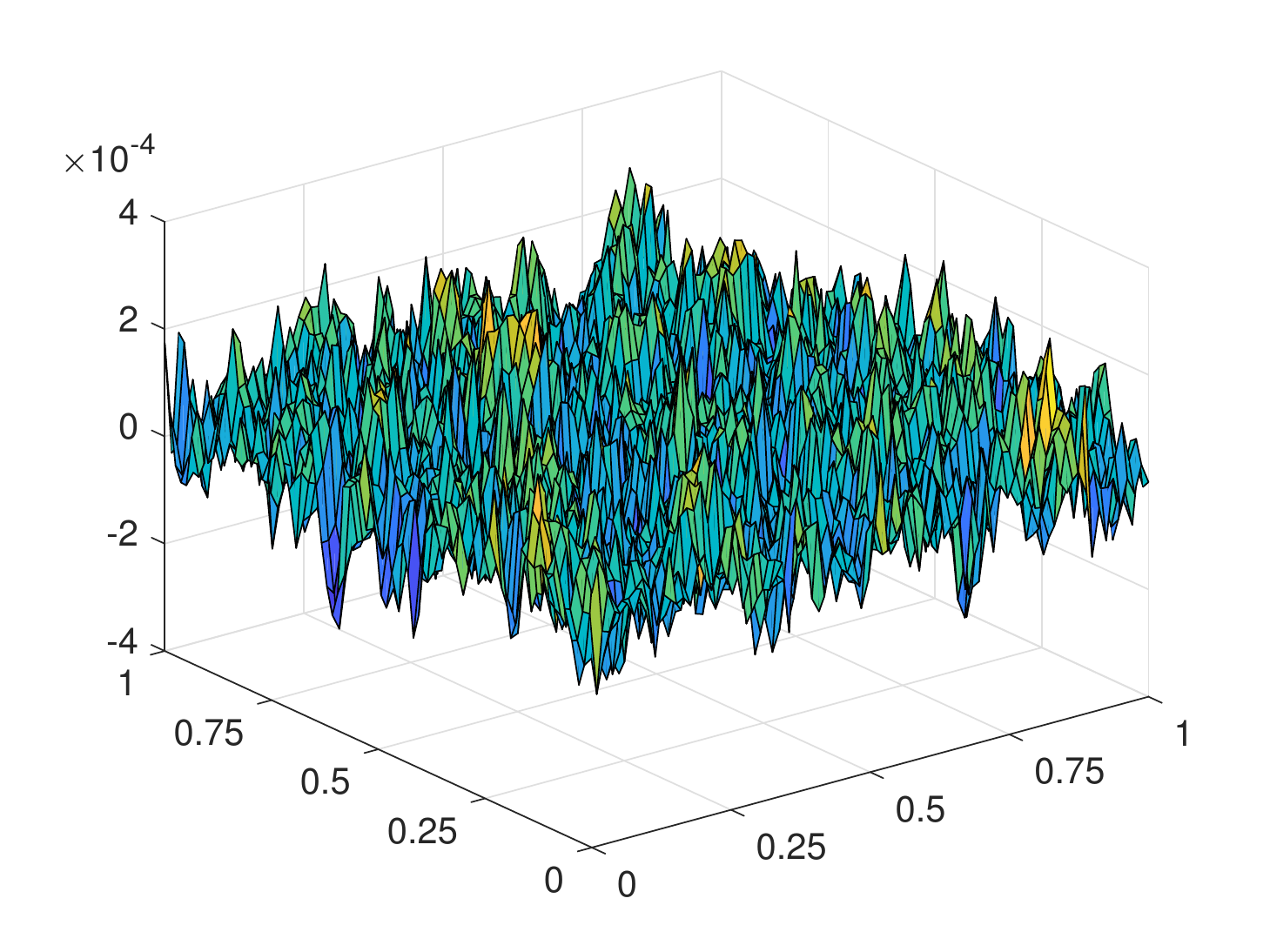}
  }
  \caption{\label{fig:se2d} A sample from the test set for the potential $\eta$, the source term
    $f$, the prediction $u^{\NN}$ with $\alpha=11$ and $\K=5$, and its corresponding error for
    Schr\"odinger form in the 2D case.}
\end{figure}

For the 2D case, the domain $\Omega=[0,1]^2$ is discretized with a $80\times 80$ uniform Cartesian
mesh. The potential $\eta(x)$ is generated by (1) sampling independently from $\cN(0,1)$ on a
uniform mesh with $10\times 10$ points, (2) then interpolating to $80\times80$ points via a Fourier
interpolation, and (3) point-wise exponentiating followed by appropriate scaling. The source term is
sampled point-wisely from a standard Gaussian distribution. When trained with $\alpha = 11$ and $\K=5$, the NN achieves a test error of
\RRev{$2.2\times 10^{-2}$ and an operator error of $4.2\times 10^{-3}$} with only $9.3\times10^5$
parameters, as reported in \cref{tab:se2d}. \RRev{The operator error $\operror$ estimate is computed by averaging the error among 10 distinct samples of the inverse operator
  $G_\eta$.} The values of $\eta$ and $f$ of a
representative sample are displayed in \cref{fig:se2d}, along with the NN prediction and the error.

\subsection{Divergence form}\label{sec:div}
The same NN architecture is applied to the Green's functions of the divergence form
\begin{equation}
  \begin{aligned}
    & -\grad\cdot(\eta(x)\grad u(x))=f(x),\quad x\in[0,1]^d,\\ & \int_{[0,1]^d}u(x)\dd x=0,
  \end{aligned}
\end{equation}
with $\eta(x)\geq\eta_0>0$ along with the periodic boundary condition. Following the notations of
\cref{sec:intro},
\begin{equation}
  L_\eta = -\grad\cdot \diag(\eta) \grad, \qquad G_\eta=L_\eta^{-1} 
  \text{ in the constraint of } \int_{[0,1]^d}u(x)\dd x=0.
\end{equation}
When $\eta(x)$ is close to a fixed $\eta_0>0$, the operator can be decomposed as
\begin{equation}
  L_0=-\eta_0\lapl, \quad E_{\eta}=\grad\cdot\diag(\eta(x)-\eta_0)\grad, \quad L_\eta = L_0-E_{\eta}.
\end{equation}
Since the operator $E_{\eta}$ is linearly dependent on $\eta$, it is easy to check that the
discussion for the Schr\"odinger form holds for the divergence form case as well.

\begin{figure}[h!]
    \centering
    \subfloat[\label{fig:Pois1d_input} $a$]{
    \includegraphics[width=0.4\textwidth, trim=.55cm 1.15cm 0.75cm 0.95cm, clip]{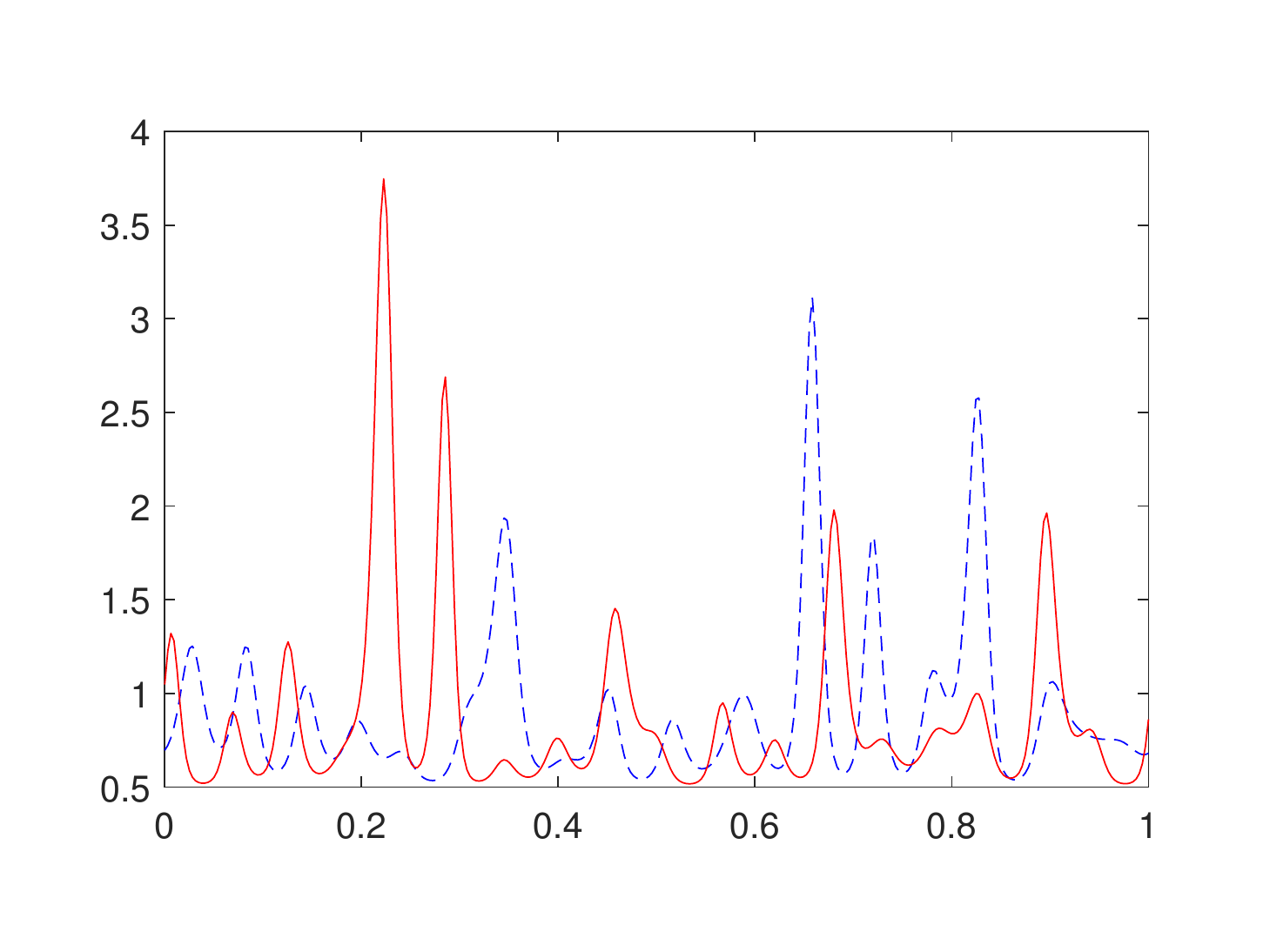}
    }
    \subfloat[\label{fig:Pois1d_input_f} $f$]{
    \includegraphics[width=0.4\textwidth, trim=.55cm 1.15cm 0.75cm 0.95cm, clip]{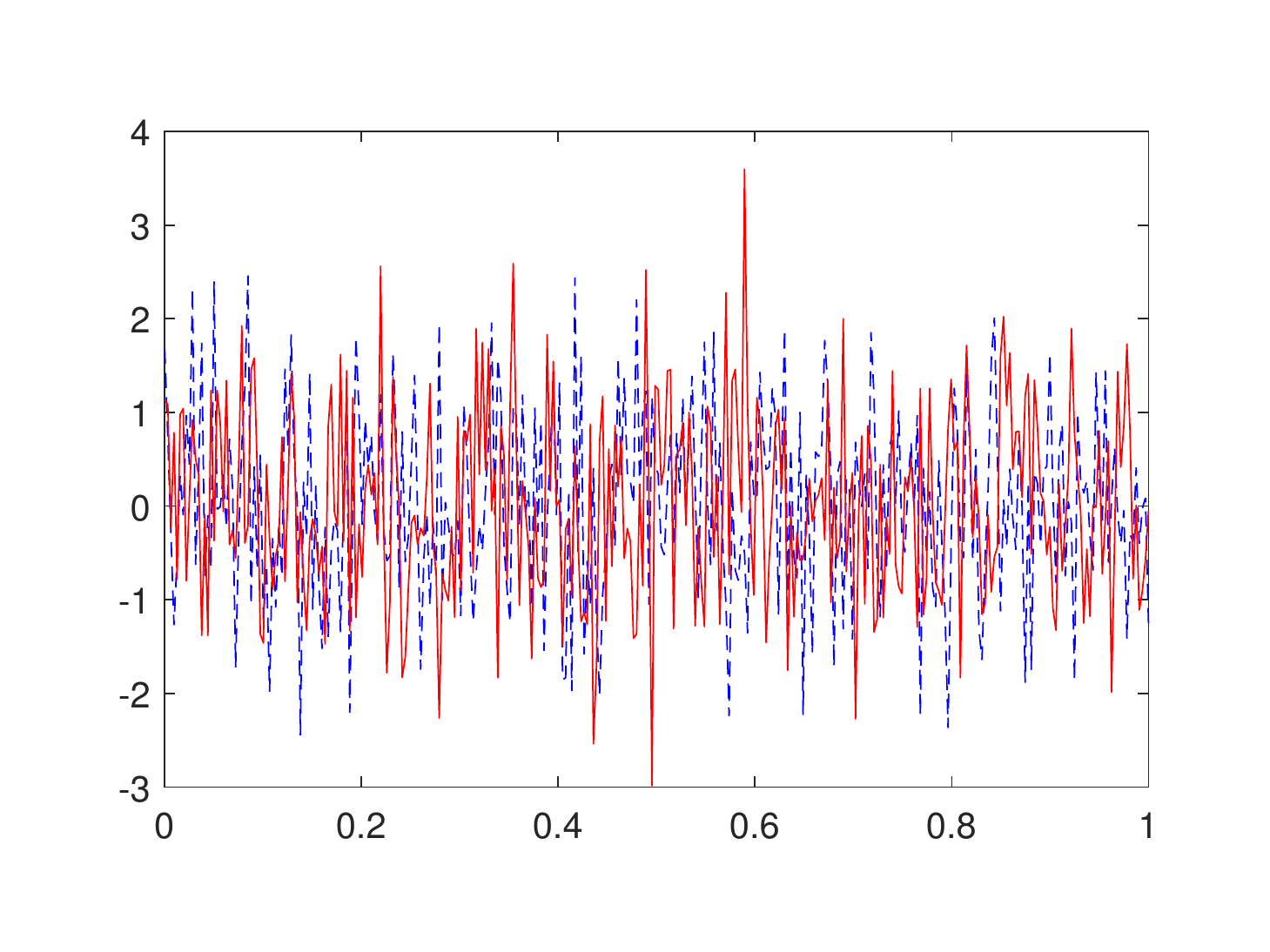}
    }\\
    \subfloat[\label{fig:Pois1d__out} $u^{\NN}$]{
    \includegraphics[width=0.4\textwidth, trim=.55cm 1.15cm 0.75cm 0.95cm,  clip]{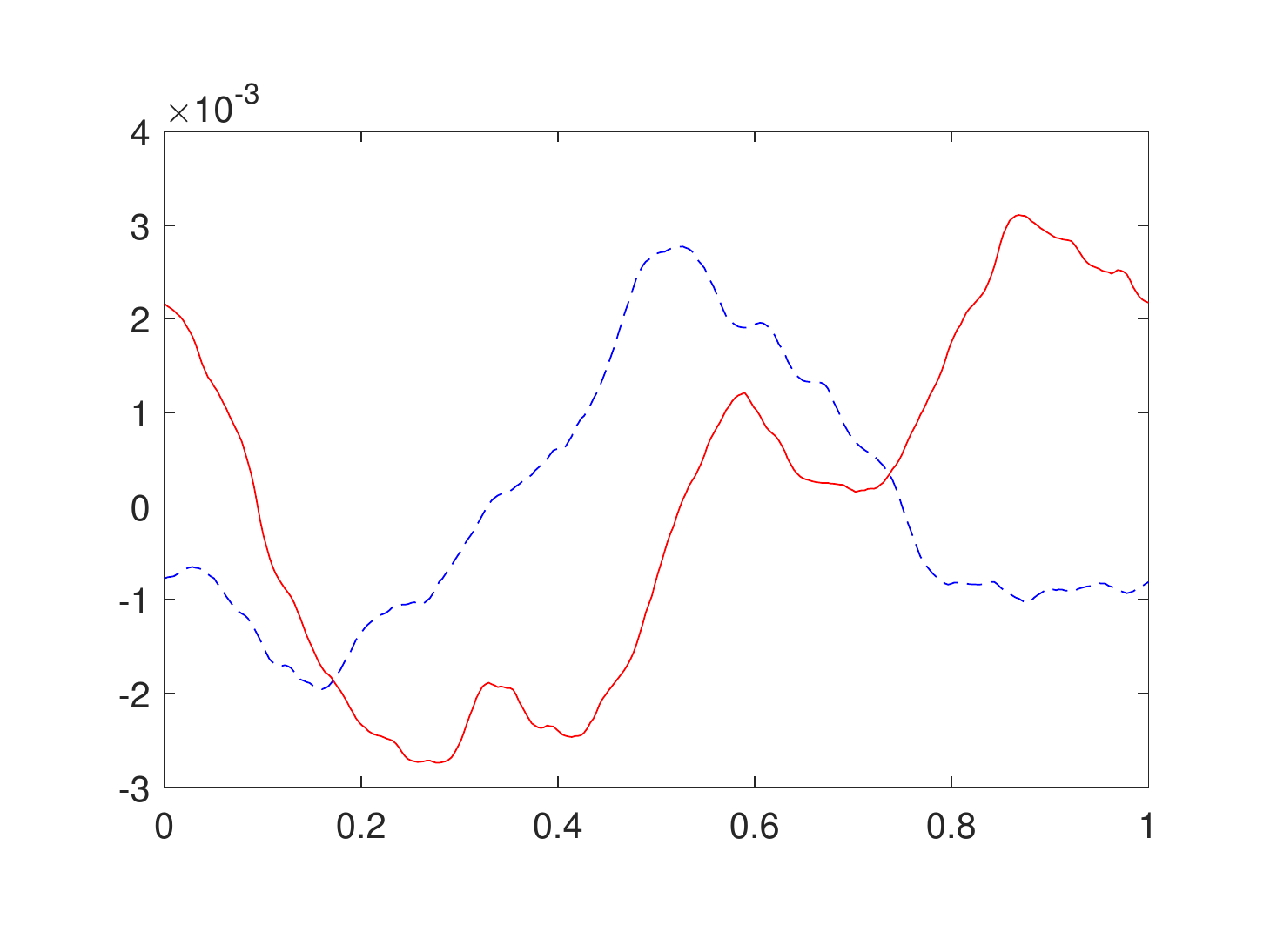}
    }
    \subfloat[\label{fig:Pois1d_err} $u^{\NN}-u$]{
    \includegraphics[width=0.4\textwidth, trim=.55cm 1.15cm 0.75cm 0.95cm, clip]{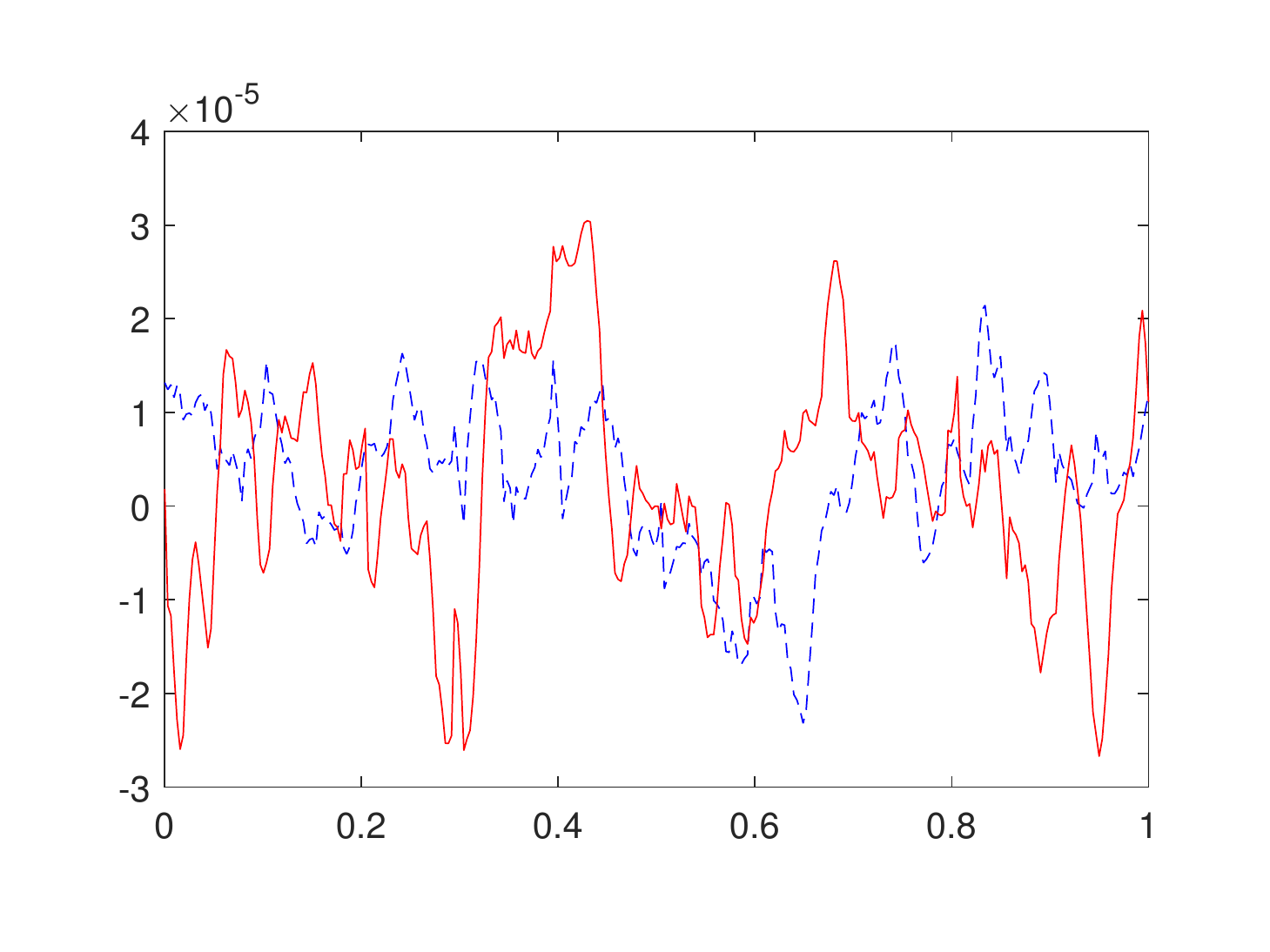}
    }
    \caption{\label{fig:pois1d} Two samples (one in red line and the other in blue dashed line) from
    the test set for the distinct coefficients $\eta$, the source terms $f$, the predictions
    $u^{\NN}$ with $\alpha=9$ and $\K=5$, and their corresponding error for the divergence form in the
    1D case.}
\end{figure}

\paragraph{Numerical results.}
The parameter field $\eta(x)$ is generated in a way similar to the potential of the Schr\"odinger
form, with the difference that the scaling factor is set to $1/5$ and an additive term of $0.5$ is
applied point-wise to avoid the ill-conditioning of $G_\eta$. The numerical results for different
choices of $\alpha$ (channel number) and $\K$ (layer number) are summarized in
\cref{tab:Pois1d}. For example, a test error of \RRev{$6.9\times10^{-3}$} is achieved at $\alpha=9$ and
$\K=5$ with \RRev{$9.6\times10^4$} parameters. Two random samples from the test data are illustrated in
\cref{fig:pois1d}.

\begin{table}[h!]
    \centering
    \begin{tabular}{ccccc}
      \hline
      $\alpha$ & $\K$ & $\Nparams$ & $\trainerror$ & $\testerror$ \\ \hline\hline
      7 & 5 &   58717  &  \RRev{7.27e-3}  &   \RRev{7.76e-3} 	 \\ \hline
      7 & 7 &   74089  &  \RRev{7.46e-3}   &     \RRev{8.29e-3}      \\ \hline
       9 & 5 &   96625 & \RRev{6.05e-3}   &    \RRev{6.88e-3}   	 \\ \hline
       9 & 7 &  122005  &  \RRev{6.83e-3}     &   \RRev{8.21e-3}    \\ \hline
    \end{tabular}
    \caption{\label{tab:Pois1d} Relative error in approximating the solution of divergence form in the 1D case.}
\end{table}

\section{Radiative transfer equation with isotropic scattering}\label{sec:rte}

The radiative transfer equation (RTE) is a fundamental model for describing particle propagation,
with applications in many fields, such as neutron transport in reactor physics
\cite{pomraning1973equations}, light transport in atmospheric radiative transfer
\cite{marshak20053d}, heat transfer \cite{koch2004evaluation}, and optical imaging
\cite{klose2002optical}. The steady-state RTE in the homogeneous scattering regime is
\begin{equation}\label{eq:RTE}
  \begin{aligned}
    v\cdot\nabla_{x}\varphi(x,v)+(\eta(x)+\eta_a(x))\varphi(x,v)
    &=\eta(x) u(x)+f(x),\quad \text{ in } \Omega\times \bbS^{d-1},\quad \Omega\subset\bbR^d,\\
    \varphi(x,v) &= 0,\quad \text{ on } \{(x,v)\in\partial\Omega\times \bbS^{d-1}: n(x)\cdot
    v<0\},\\
    u(x) &= \frac{1}{4\pi}\int_{\bbS^{d-1}}\varphi(x,v)\dd v,
  \end{aligned}
\end{equation}
where $\varphi(x,v)$ denotes the photon flux that depends on both space $x$ and angle $v$, $f(x)$ is
the light source, $\eta(x)$ is the scattering coefficient, and $\eta_a(x)$ is the physical
absorption coefficient. In many applications, it is reasonable to assume $\eta_a(x)$ to be
constant. Below, we focus on the most challenging case $\eta_a(x) \equiv 0$.

The numerical solution to the RTE has been extensively studied using the Monte Carlo methods and
various discretization schemes for the differential-integral formulation \cref{eq:RTE} of
RTE. However, these approaches often suffer from the high-dimensionality and non-smoothness of the
photon-flux $\varphi(x,v)$. The recent numerical work in \cite{LionsVol6,fan2018fast,ren2016fast}
follows the integral formulation by eliminating $\varphi(x,v)$ from the equation and keeping only $u(x)$ as unknown:
\begin{equation}\label{eq:RTE_IE}
    \left( I - K_\eta \eta \right) u = K_\eta f,
\end{equation}
where the operator $K_\eta$ is defined as
\begin{equation}
  K_\eta f = \int_{y\in\Omega}k_\eta(x,y)f(y)\dd y,\quad
  k_\eta (x,y) = \frac{\exp\left( -|x-y|\int_0^1\eta(x-s(x-y))\dd s \right)}{4\pi |x-y|^{d-1}}.
\end{equation}
The parameterized Green's function operator for the steady-state RTE is then 
\begin{equation}
  \label{eq:Gmu}
  \Gmu = \left( I -K_\eta \eta \right)^{-1}K_\eta.
\end{equation}
Since $K_\eta$ is a dense operator, forming $\Gmu$ following \cref{eq:Gmu} is often computationally
expensive. Instead, the meta-learning approach developed above allows for approximating the map from
$\eta$ to $\Gmu$ directly.

\cref{sec:sch} argues that the map $\eta\to C_{\eta}\sps{\ell}$ for the translation invariant
operator can be represented by a convolutional NN.  A key observation for the current setting is
that the integral equation \cref{eq:RTE_IE} can be extended to the whole domain by padding $f$ and
$\eta$ with zero. As a result, the map from $\eta$ to $G_{\eta}$ can be represented by a
convolutional NN {\em with zero padding}.

\paragraph{Numerical results.}

\begin{figure}[h!]
    \centering
    \subfloat[\label{fig:rte1d_input} $\eta$]{
    \includegraphics[width=0.4\textwidth, trim=.95cm 1.15cm 0.75cm 1.15cm, clip]{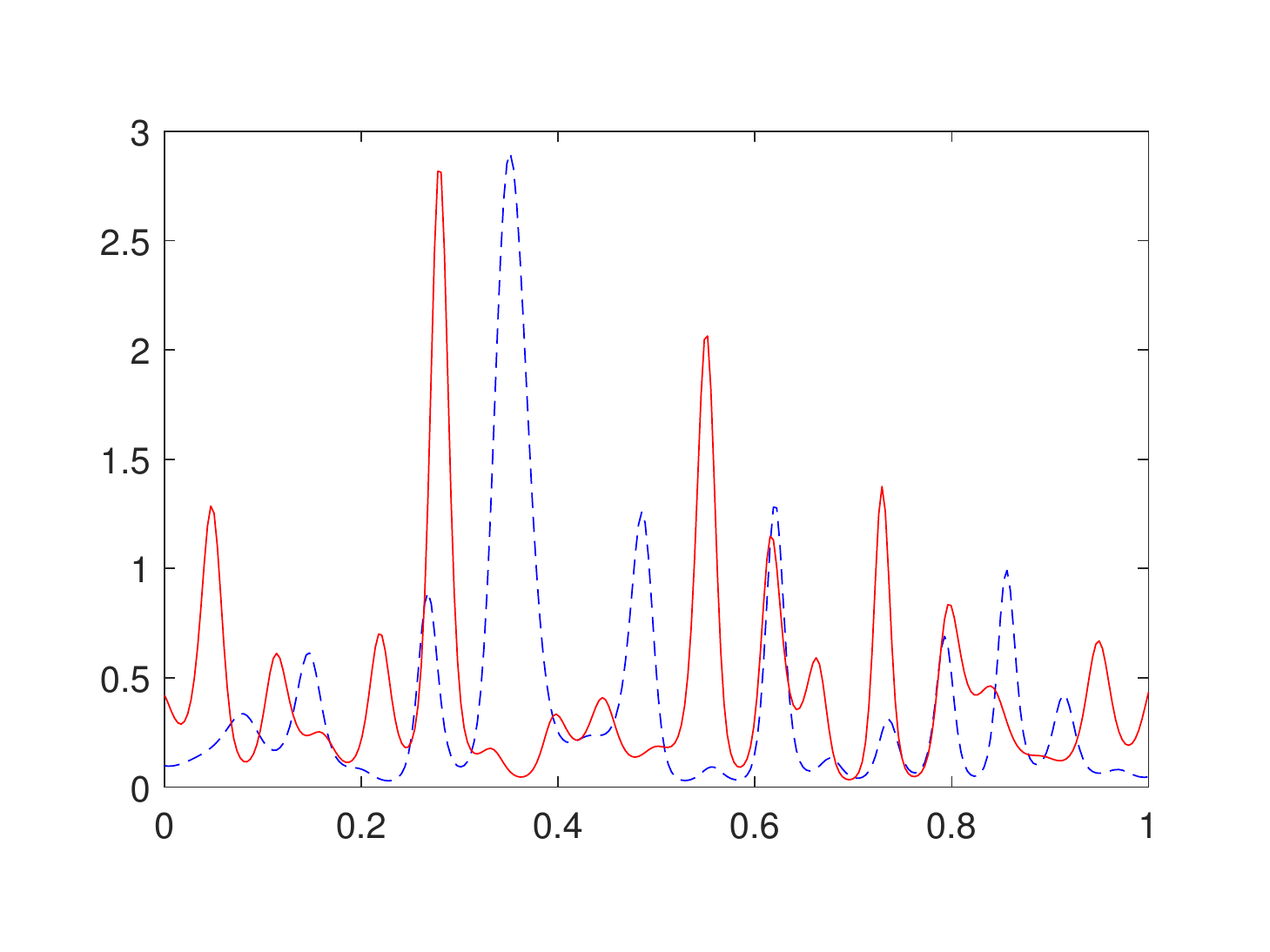}
    }
    \subfloat[\label{fig:rte1d_input_f} $f$]{
    \includegraphics[width=0.4\textwidth, trim=.95cm 1.15cm 0.75cm 1.15cm, clip]{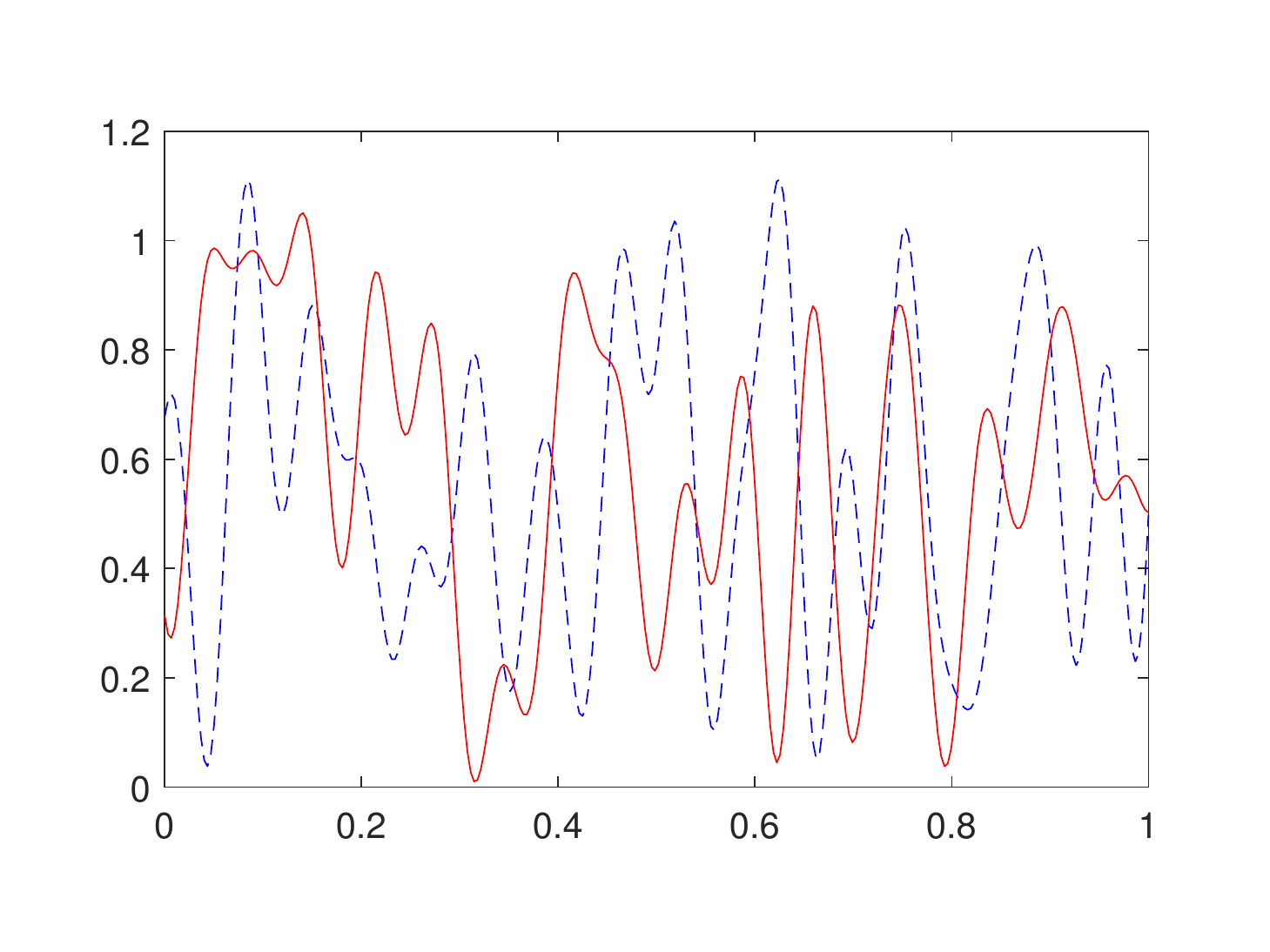}
    }\\
    \subfloat[\label{fig:rte1d__out} $u^{\NN}$]{
    \includegraphics[width=0.4\textwidth, trim=.95cm 1.15cm 0.75cm 1.15cm,  clip]{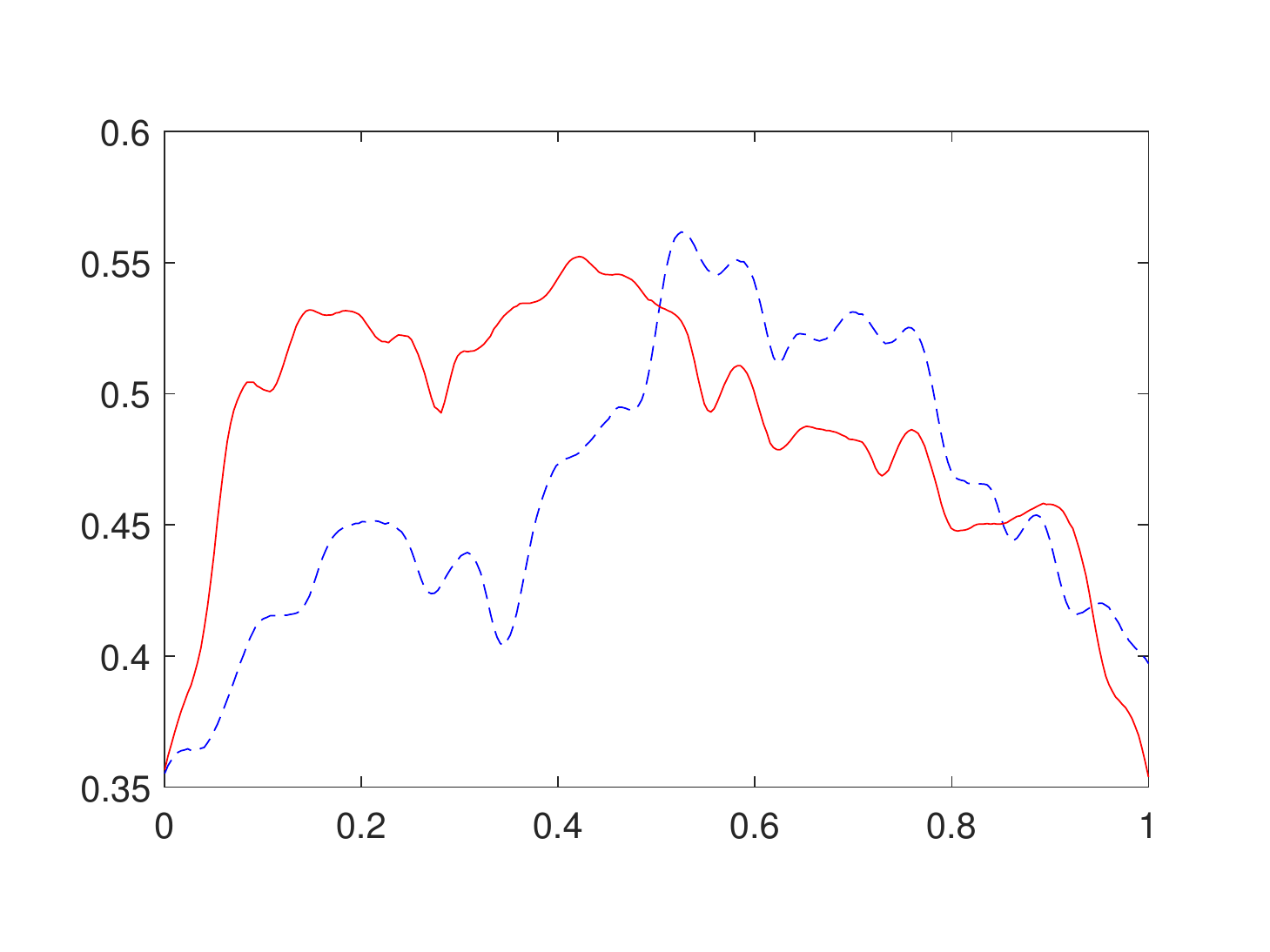}
    }
    \subfloat[\label{fig:rte1d_err} $u^{\NN}-u$]{
    \includegraphics[width=0.4\textwidth, trim=.95cm 1.15cm 0.75cm 0.95cm, clip]{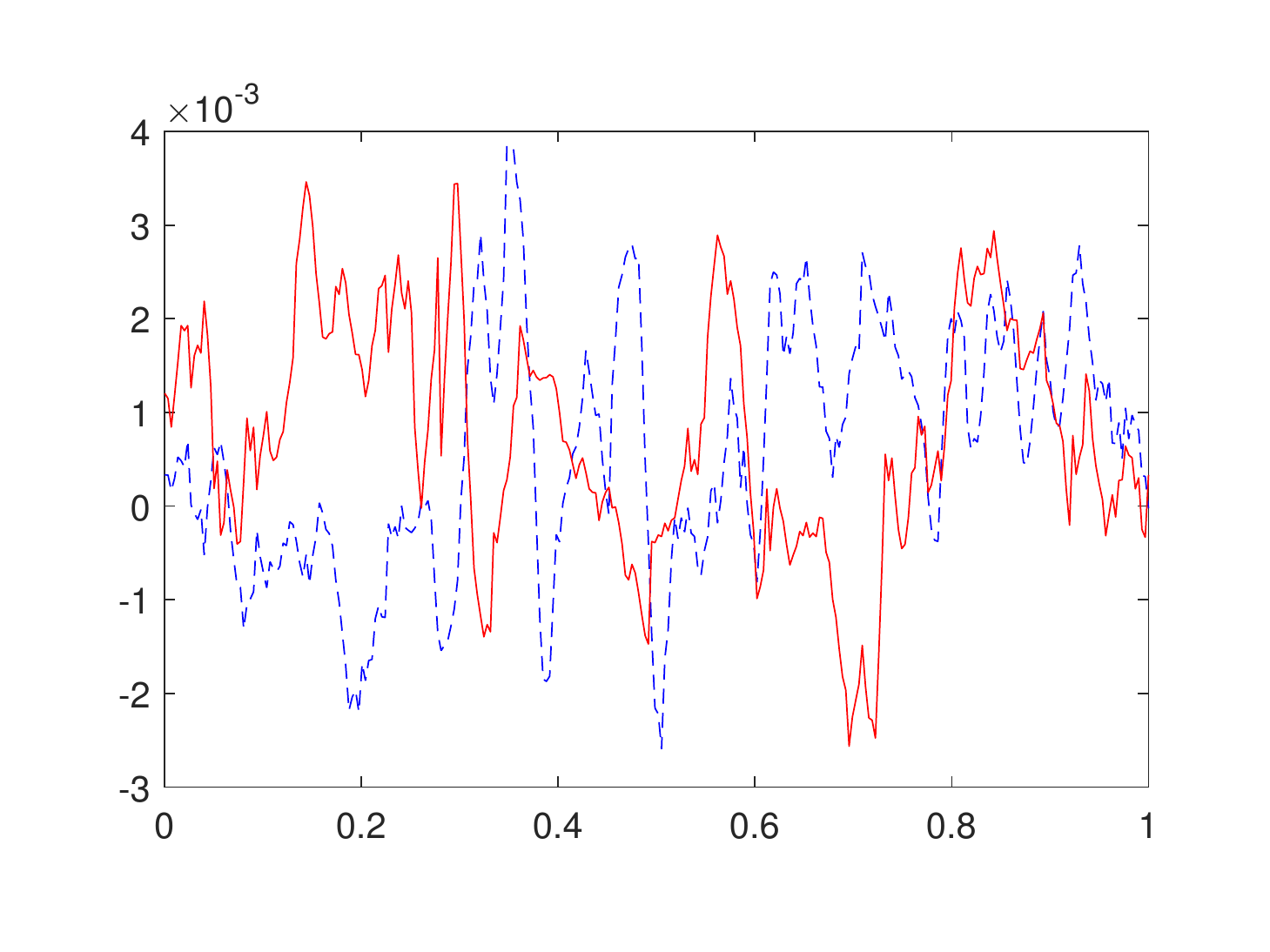}
    }
    \caption{\label{fig:rte1d} Two samples (one in red line and the other in blue dashed line) from
    the test set with the scattering coefficients $\eta$, the source terms $f$, the predictions
    $u^{\NN}$ with $\alpha=5$ and $\K=5$, and their corresponding error in the 1D RTE.}
\end{figure}

The first test is concerned with the one-dimensional slab geometry, where the parameter $\eta$ varies
only in the $x_1$ direction (i.e., constant in the $x_2$ and $x_3$ directions). For this geometry,
the integral equation \cref{eq:RTE_IE} reduces to
\begin{equation}\label{eq:RTE_IE1d}
  (I-K^{(1)}_\eta \eta)u(x)=K^{(1)}_\eta f(x),\quad
  G_\eta = (I-K^{(1)}_\eta \eta)^{-1} K^{(1)}_\eta,
\end{equation}
where $x$ stands for only $x_1$ and the operator $K_\eta^{(1)}$ is defined as
\begin{equation}
  \begin{aligned}
    K^{(1)}_\eta f(x) &= \int_{y\in \Omega} k^{(1)}_\eta (x,y)f(x)\dd y,\\
    k^{(1)}_\eta(x,y) &= \frac{1}{2}\mathrm{Ei}\left(-|x-y|\int_{0}^1\eta(x-s(x-y))\dd s\right),
    \quad \mathrm{Ei}(x)=\int_{-\infty}^{x}\frac{e^t}{t}\dd t,
  \end{aligned}
\end{equation}
with the domain $\Omega=[0,1]$. In the implementation, $[-x_0,1+x_0]$ is discretized by a uniform
Cartesian mesh with $N=320$ points, where $x_0>0$ is selected such that there are $300$ points in
$\Omega$.  The scattering coefficient $\eta$ is generated in the same way for $\eta(x)$ in
\cref{sec:sch} followed by appropriate rescaling. The source term $f(x)$, positive due to physical
considerations, is generated by sampling independently from $\cU(0,1)$ instead of $\cN(0,1)$ and interpolated via Fourier
interpolation. The
values of $\eta$ and $f$ outside of $\Omega$ are set to be $0$. The results for different values of
$\alpha$ (channel number) and $\K$ (layer number) are summarized in \cref{tab:rte1d}.  A test error
of $2.9 \times 10^{-3}$ is achieved with as few as $3.4\times 10^4$ parameters with
$\alpha=\K=5$. Two representative examples from the test set are shown in \cref{fig:rte1d}.

\begin{table}[h!]
    \centering
    \begin{tabular}{ccccc}
      \hline
      $\alpha$ & $K$ & $\Nparams$ & $\trainerror$ & $\testerror$ \\ \hline\hline
       5 & 5  &  34131  & 2.48e-3 & 2.93e-3 	\\ \hline
      5 & 7  &   41991 &  2.46e-3 & 3.01e-3   \\ \hline
      7 & 5 &  66403  &  1.92e-3 & 2.45e-3  \\ \hline
      7 & 7 &  81775  &  2.05e-3 & 2.36e-3 \\ \hline
                \end{tabular}
    \caption{\label{tab:rte1d} Relative error in approximating the solution to the 1D RTE.}
\end{table}

\begin{figure}[h!]
    \centering
    \subfloat[\label{fig:rte2d_input} $V$]{
    \includegraphics[width=0.4\textwidth,trim=0.5cm 0.5cm 0.9cm 0.7cm,  clip]{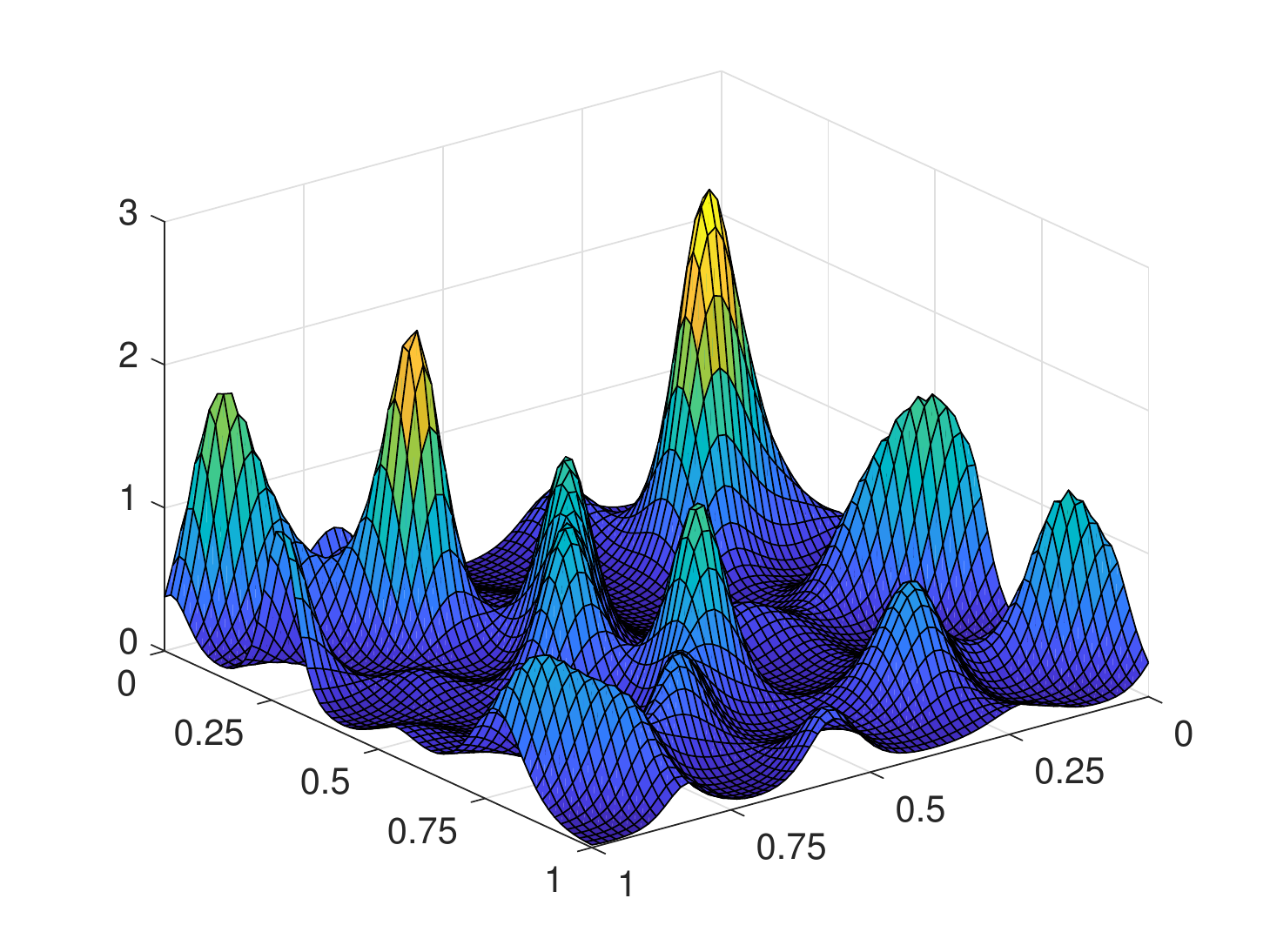}
    }
    \subfloat[\label{fig:rte2d_input_f} $f$]{
    \includegraphics[width=0.4\textwidth,trim=0.5cm 0.5cm 0.9cm 0.7cm,  clip]{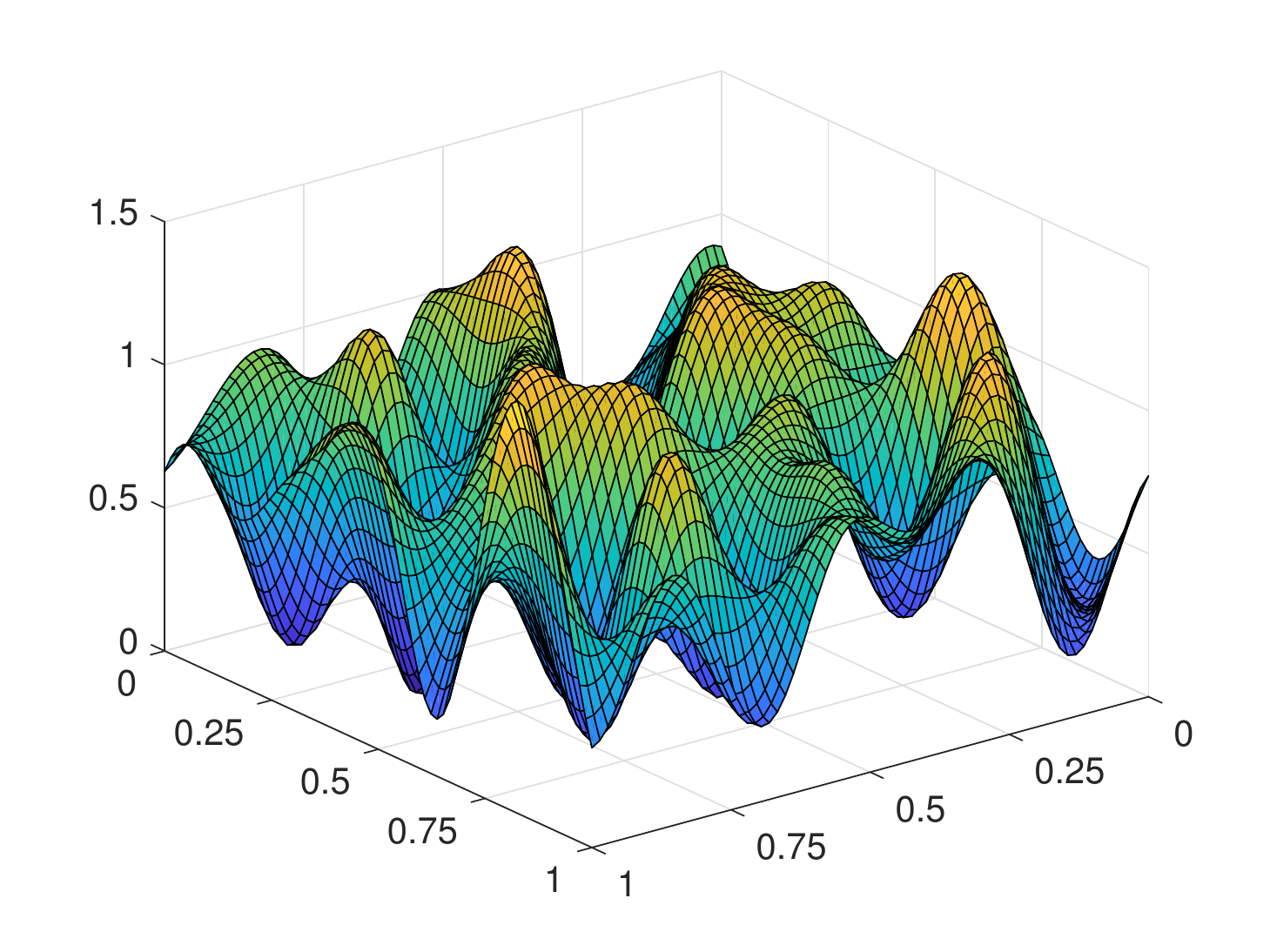}
    }\\
    \subfloat[\label{fig:rte2d__out} $u^{\NN}$]{
    \includegraphics[width=0.4\textwidth,trim=0.5cm 0.5cm 0.9cm 0.7cm,  clip]{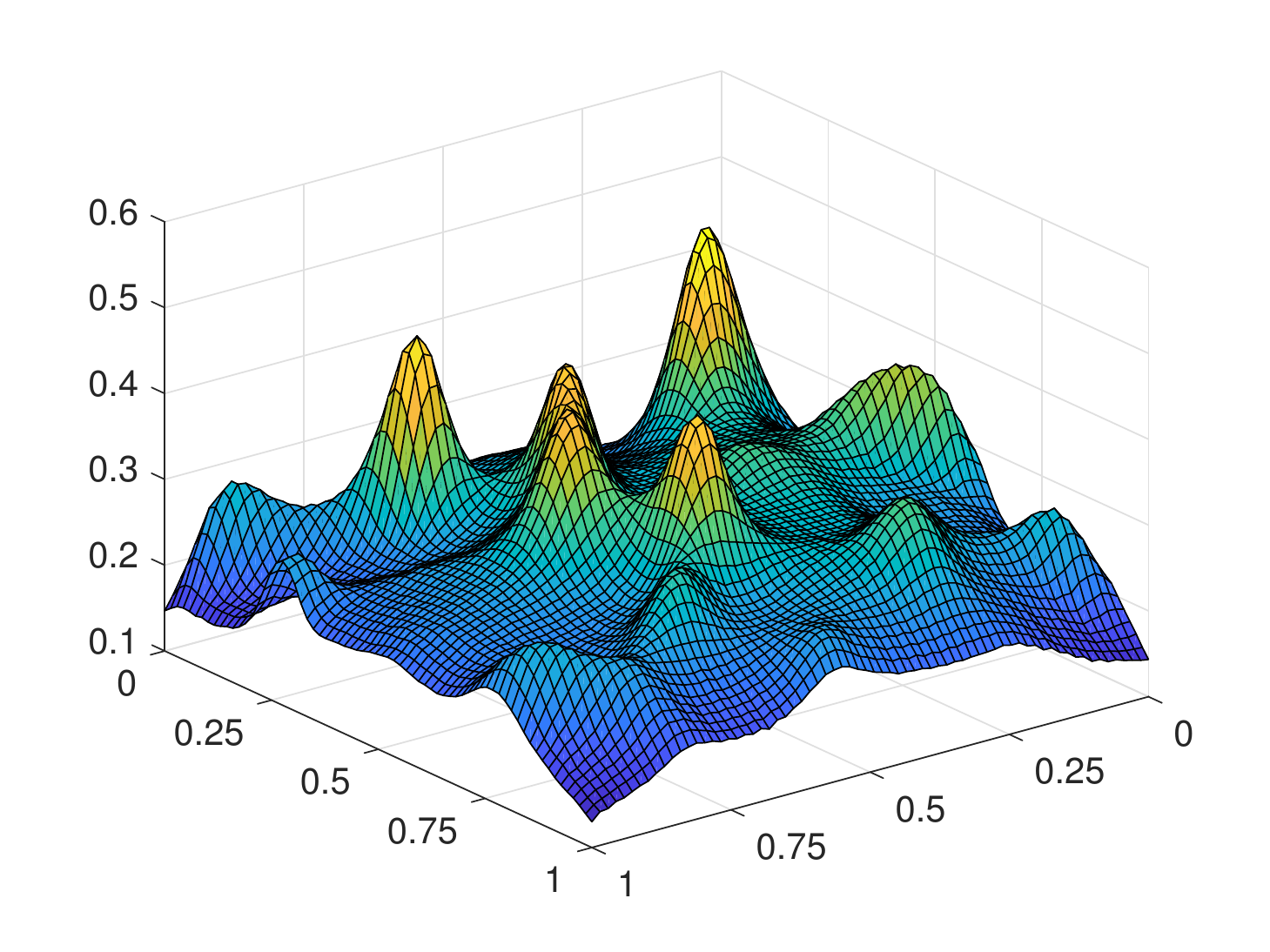}
    }
    \subfloat[\label{fig:rte2d_err} $u^{\NN}-u$]{
    \includegraphics[width=0.4\textwidth,trim=0.5cm 0.5cm 0.9cm 0.7cm,  clip]{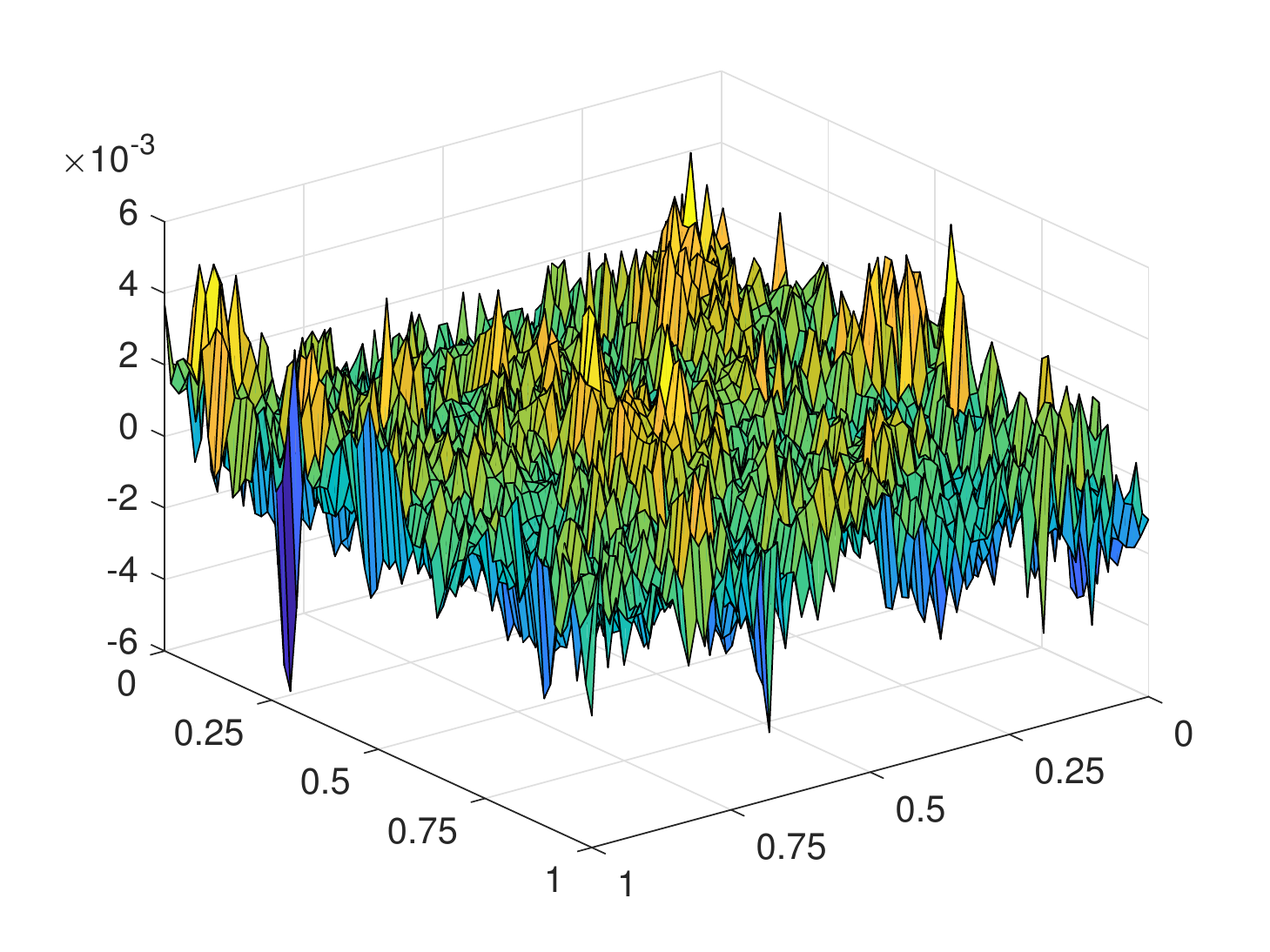}
    }
    \caption{\label{fig:rte2d} A sample from the test set for the scattering coefficient $\eta$, the
    source term $f$, the prediction $u^{\NN}$ with $\alpha=11$ and $\K=5$, and the corresponding
    error for the 2D RTE.}
\end{figure}

The second test is concerned with the 2D RTE. The domain $\Omega = [-x_0,1+x_0]^2$ 
is discretized with a uniform Cartesian grid with $80\times 80$ points, where $x_0$ is chosen such
that there are $70\times 70$ points in $\Omega$. The scattering coefficient is generated following
the same way of $\eta$ in \cref{sec:sch} for the 2D case, followed by an appropriate rescaling.  The
source term $f(x)$ is generated by sampling independently from $\cU(0,1)$ instead of $\cN(0,1)$.
The values of $\eta$ and $f$ outside of $\Omega$ are set to be $0$.  Results reported in
\cref{tab:rte2d} show that by setting $\alpha=11$ and $\K=5$, the NN can achieve a test error of
$4.4\times 10^{-3}$ with as few as $1.3\times 10^6$ parameters. A representative sample from the
test set is illustrated in \cref{fig:rte2d}.

\begin{table}[h!]
    \centering
    \begin{tabular}{ccccc}
      \hline
      $\alpha$ & $K$ & $\Nparams$ & $\trainerror$ & $\testerror$ \\ \hline\hline
      11 & 5  & 1287903 &  4.39e-3 & 4.39e-3	\\
      15 & 5  &  1663831 & 3.55e-3 & 3.55e-3 	\\ \hline
    \end{tabular}
    \caption{\label{tab:rte2d} Relative error in approximating the solution to the 2D RTE.}
\end{table}

\section{Conclusions}\label{sec:conc}

This paper presented a meta-learning approach for learning the map from the equation parameter $\eta$
to the pseudo-differential solution operator $G_\eta$. Motivated by the nonstandard wavelet form
\cite{bcr}, the pseudo-differential operator is compressed to a collection of vectors. The nonlinear
map from the parameter to this collection of vectors and the wavelet transform are learned
hand-in-hand in the meta-learning approach. Numerical studies are carried out for the Green's
functions of elliptic PDEs as well as the radiative transfer equation.

This approach can be extended in several directions. First, this paper is only concerned with linear
operators $G_\eta$. This work can be readily extended to nonlinear operators if a simple compressed
representation (such as the collection of vectors used here) can be identified. Second, the \ConvNet
module for the map $\eta\to C_{\eta}\sps{\ell}$ can be replaced with the recently proposed
multiscale NNs \cite{bcrnet, fan2018mnn2, fan2018mnn}, which are more effective for certain
global-scale convolutions.

\section*{Acknowledgments}
The work of Y.F. and L.Y. is partially supported by the U.S. Department of Energy, Office of
Science, Office of Advanced Scientific Computing Research, Scientific Discovery through Advanced
Computing (SciDAC) program.  The work of J.F. is partially supported by Stanford Graduate Fellowship
in Science \& Engineering and by ``la Caixa'' Fellowship, sponsored by the ``la Caixa'' Banking
Foundation of Spain under  Fellowship LCF/BQ/AA16/11580045. The work of L.Y. is also partially
supported by the National Science Foundation under award DMS-1818449.
This work is also supported by the GCP Research Credits Program from Google and AWS Cloud Credits
for Research program from Amazon.

\bibliographystyle{abbrv}
\bibliography{references}
\end{document}